\documentclass[10pt]{article}
\usepackage[ansinew]{inputenc}
\usepackage{amsmath}
\usepackage[francais]{babel}
\usepackage[pdftex]{graphicx}
\title{Cohomologie $L^p$ et pincement
\footnote{Mots cl{\'e} : Cohomologie $L^{p}$, courbure n{\'e}gative, espace
homog{\`e}ne, espace de Besov. Keywords : $L^{p}$-cohomology, negative curvature, homogeneous space, Besov space.
Mathematics Subject Classification :
43A15, %$L^{p}$ and other function spaces on groups
43A80, %Analysis on non compact Lie groups
46E35, %Sobolev spaces and other spaces of smooth functions
53C20, %Global Riemannian geometry, including pinching
53C30, %Homogeneous manifolds
58A14. %Hodge theory
}}
\author{P. Pansu}

\newtheorem{theo}{Th{\'e}or{\`e}me}
\newtheorem{lemme}{Lemme}
\newtheorem{prop}[lemme]{Proposition}
\newtheorem{cor}[lemme]{Corollaire}
\newtheorem{defi}[lemme]{D{\'e}finition}
\newtheorem{rem}[lemme]{Remarque}
\newtheorem{exemple}[lemme]{Exemple}
\def\preuve{\par\medskip\noindent {\bf Preuve.}{\hskip1em}}
\def\remarque{\par\medskip\noindent {\bf Remarque.}{\hskip1em}}

\def\qed{~q.e.d.}

\def\R{{\bf R}}
\def\C{{\bf C}}
\def\Z{{\bf Z}}
\def\ltimes#1{\times_{\alpha}}
\def\n#1{{\parallel #1 \parallel}}
\def\l{\lambda}

\def\L{\Lambda}

\def\H{{\cal H}}

\def\e{\epsilon}

\def\dim{{\rm dim}\,}
\def\ker{{\rm ker}\,}
\def\im{{\rm im}\,}
\def\tr{{\rm tr}\,}

\def\ddt{{{\partial}\over{\partial t}}}
\def\con{\hbox{const.}}
\def\d{\displaystyle}
\def\inv#1{\frac{1}{#1}}
\def\op#1{\Omega^{#1,p}}
\def\opr#1{\Omega^{#1,p'}}
\def\hp#1{H^{#1,p}}
\def\hop#1{H_{o}^{#1,p}}
\def\tp#1{T^{#1,p}}
\def\toop#1{T_{o}^{#1,p}}
\def\rp#1{R^{#1,p}}

\def\lp{_{L^{p}}}
\def\et{\quad\hbox{et}\quad}
\def\si{\quad\hbox{si}\quad}
\def\sinon{\quad\hbox{sinon}}
\def\ali{\medskip}

\def\dplus{d_{+}}
\def\dmoins{d_{-}}
\def\dzero{d_{0}}
\def\timez{\times\hskip-2pt\vrule height 1ex width .055ex depth -.0ex}

\begin{document}
\maketitle
\begin{quote}
{\small
RESUME. On donne un crit{\`e}re optimal d'annulation de la torsion en cohomologie $L^p$ pour les vari{\'e}t{\'e}s riemanniennes {\`a} courbure sectionnelle n{\'e}gative pinc{\'e}e. Il en résulte que certains espaces homogènes à courbure négative ne sont pas quasiisométriques à des variétés plus pincées qu'eux.

ABSTRACT. A sharp vanishing theorem for the $L^p$ cohomology torsion of Riemannian manifolds with pinched negative curvature is given. It follows that certain negatively curved homogeneous spaces cannot be quasiisometric to better pinched manifolds.}
\end{quote}

\section{Introduction}

\subsection{Motivation : un probl\`eme de pincement}

D'un th\'eor\`eme de M. Berger et W. Klingenberg, 
\cite{Berger}, il r\'esulte que si $V$ est un espace
sym\'etrique de rang un de type compact \`a courbure non
constante (i.e. un espace projectif complexe $\C P^{m}$, $m\geq 2$, un espace projectif quaternionien ${\bf H} P^{m}$, $m\geq 2$, ou le plan projectif des octaves de Cayley ${\bf Ca}P^{2}$), $V$ n'admet pas de m\'etrique \`a courbure comprise entre $\delta$ et $1$ si $\delta>\inv4$. 

\bigskip

On se pose un probl{\`e}me analogue en courbure n{\'e}gative. Si $-1\leq\delta<0$, on dit qu'une variété riemannienne est {\em $\delta$-pinc\'ee} s'il existe $a>0$ tel que sa courbure soit comprise entre $-a$ et $\delta a$. 

Par exemple, l'espace hyperbolique réel est $-1$-pincé. Les espaces
sym\'etriques de rang un de type non compact \`a courbure non
constante sont $-\inv4$-pinc{\'e}s. Il s'agit des espaces hyperboliques complexes $\C H^{m}$, $m\geq 2$, des espaces hyperboliques quaternioniens ${\bf H} H^{m}$, $m\geq 2$, et du plan hyperbolique des octaves de Cayley ${\bf Ca}H^{2}$.

\bigskip

Le problème du pincement optimal consiste à déterminer quel est le meilleur pincement possible pour une métrique sur une variété donnée. Pour les variétés simplement connexes (et donc difféomorphes à l'espace hyperbolique réel), il convient de se restreindre {\`a} des m{\'e}triques comparables {\`a} une m{\'e}trique de référence, par exemple, qui lui sont quasiisom{\'e}triques. On rappelle que deux vari\'et\'es riemanniennes $M$ et $N$ sont dites {\sl quasi\-i\-so\-m\'e\-tri\-ques} s'il existe une application $f:M\to N$ et des constantes $C$ et $L$ telles que l'image de $f$ soit $C$-dense dans $N$ et pour tous points $x$, $y\in M$, 
$$
-C+\inv{L}\leq d(f(x),f(y))\leq Ld(x,y)+C.
$$

\bigskip

{\bf Question}. {\em Soit $M$ une variété riemannienne $\delta$-pinc\'ee. Existe-t'il une vari\'et\'e riemannienne $N$ quasiisom\'etrique \`a $M$ et $\delta'$-pinc\'ee avec $\delta'<\delta$ ?}

\bigskip

Dans cet article, on détermine le pincement optimal pour des familles
d'espaces homog\`enes riemanniens. Voici un exemple. Soit $G_{2,4,-\inv4}$ le produit semi-direct de $\R^3$ par $\R$ d\'efini par le groupe \`a un
param\`etre d'automorphismes de $\R^3$ engendr\'e par la matrice  
$$\begin{pmatrix}1 & 0 & 0 \\ 0 & 1 & 0\\ 0 & 0 & 2\end{pmatrix}.$$ 
La m\'etrique riemannienne qui en coordonn\'ees exponentielles $t$ (sur le facteur $\R$), $x$, $y$ et $z$ (sur le facteur $\R^3$) s'\'ecrit
$$
ds^{2}=dt^{2}+e^{2t}(dx^{2}+dy^{2})+e^{4t}dz^{2}
$$
est invariante \`a gauche. On v\'erifie ais\'ement (voir par exemple \cite{Heintze}) que cette m\'etrique est $-\inv4$-pinc\'ee.

\begin{theo}
\label{1/4}
Soit $\delta<-\inv4$. Aucune vari\'et\'e riemannienne
$\delta$-pinc\'ee n'est quasiisom\'etrique \`a $G_{2,4,-\inv4}$.
\end{theo}

La preuve utilise la {\sl torsion en cohomologie $L^p$}.
C'est un espace vectoriel, not\'e $\tp2(M)$, d\'efini pour $p\geq
1$. Pour une vari\'et\'e simplement connexe \`a courbure
n\'egative, le nombre 
$$
T(M)=\inf\{p>1~;~\tp2(M)\not=0\}
$$
est un invariant de quasiisom\'etrie. Un th\'eor\`eme de
comparaison (th\'e\-o\-r\`e\-me A) entra\^{\i}ne que si $\dim M=4$ et si
$M$ est $\delta$-pinc\'ee, alors $\d T(M)\geq
1+2\sqrt{-\delta}$. Un calcul direct (th\'eor\`eme B) montre
que pour le produit semi-direct $G_{2,4,-\inv4}$, la
torsion $\tp2$ est non nulle pour $2<p<4$, d'o\`u 
$$
T(G_{2,4,-\inv4})=2.
$$
La minoration du pincement s'en d\'eduit imm\'ediatement.

\subsection{Un problème ouvert}

A ma connaissance, le problème du pincement optimal pour les espaces symétriques $-\inv4$-pincés est toujours ouvert. Pourtant, le plan hyperbolique complexe $\C H^2$ est infiniment voisin de $G_{2,4,-\inv4}$. Il peut-\^etre vu comme un groupe de Lie r\'esoluble muni d'une m\'etrique
invariante \`a gauche. Ce groupe est le produit
semi-direct du groupe de Heisenberg $Heis$ par $\R$
engendr\'e par la d\'erivation de matrice 
$\d\begin{pmatrix}1 & 0 & 0 \\ 0 & 1 & 0\\ 0 & 0 & 2\end{pmatrix}$. 
Toutefois 
$$
T(\C H^{2})=4,
$$
si bien que le th\'eor\`eme de comparaison ne donne pas de
borne optimale pour le pincement des vari\'et\'es riemanniennes $N$
quasiisom\'etriques au plan hyperbolique complexe. Il y a donc une limitation essentielle dans la m\'ethode.

\ali

Le probl{\`e}me restreint o{\`u} l'on suppose que la variété inconnue $N$ rev\^et une vari{\'e}t{\'e} riemannienne compacte a {\'e}t{\'e} r{\'e}solu par M. Ville \cite{Ville} en dimension $4$, par L. Hern{\'a}ndez
\cite{Hernandez}, S.T. Yau et F. Zheng, \cite{YZ} pour les espaces hyperboliques complexes, par N. Mok, Y.T. Siu et S.K. Yeung \cite{MSY}, J. Jost et S.T. Yau \cite{JY} pour les autres espaces symétriques de rang 1.

\subsection{Cohomologie $L^{p}$}

Soit $M$ une vari{\'e}t{\'e} riemannienne. Soit $p>1$ un r{\'e}el. On note
$L^{p}\Omega^{*}(M)$
l'espace de Banach des formes diff{\'e}rentielles $L^p$ et
$\Omega^{*,p}(M)=L^{p}\cap
d^{-1}L^{p}$ l'espace des formes diff{\'e}rentielles
$L^p$ dont la diff{\'e}rentielle ext{\'e}rieure est aussi $L^p$. La cohomologie du complexe $(\Omega^{*,p}(M),d)$ s'appelle la {\em cohomologie} $L^p$ de $M$. Elle est int{\'e}ressante surtout si $M$ est non compacte. 

\ali

Par d\'efinition, la cohomologie $L^p$ est invariante par
diff\'eomorphisme bilipschitzien. Dans la classe des
vari\'et\'es simplement connexes \`a courbure n\'egative ou nulle,
c'est un invariant de quasiisom\'etrie (cf. \cite{Gromovasympt}). 

En toute g{\'e}n{\'e}ralit{\'e}, la cohomologie $L^p$ se d{\'e}compose en
cohomologie r{\'e}duite et torsion
$$ 
0\to \tp{*}\to\hp{*}\to\rp{*}\to 0,
$$
o{\`u} la {\em cohomologie r{\'e}duite} est $\rp{*}=\ker d/\overline{\im d}$ et
la {\em torsion} est $\tp{*}=\overline{\im d}/\im d$. La cohomologie
r{\'e}duite (parfois notée $\overline{H}^{k}_{(p)}$) est un espace de Banach sur lequel les isom{\'e}tries de $M$ agissent isom{\'e}triquement. La torsion est non s{\'e}par{\'e}e.

Par exemple, la cohomologie $L^p$ de la droite r\'eelle est
enti\`erement de torsion. La cohomologie $L^p$ du plan
hyperbolique est enti\`erement r\'eduite. N\'eanmoins,
cohomologie r\'eduite et torsion coexistent souvent.

\subsection{Pincement de la courbure}

En degr{\'e}s $k>1$, la cohomologie
$L^p$ est li{\'e}e de fa\c con optimale au pincement de la courbure.

\begin{stheo} {\bf A.~}
\label{A}
{\em
Soient $\delta\in]-1,0[$ un r{\'e}el, $n$ 
et $k=2,\ldots,n$ des entiers. Notons $\d{\bf
q}(n,\delta,k)=1+{{n-k-1}\over{k}}\sqrt{-\delta}$. 

Soit $M$ une vari{\'e}t{\'e} riemannienne compl{\`e}te de dimension $n$,
simplement connexe, dont la courbure sectionnelle $K$ satisfait $-1\leq K\leq\delta$. 
Alors
$$
\tp{k}(M)=0, \quad\hbox{i.e.}\quad \hp{k}(M)\quad\hbox{est s{\'e}par{\'e} pour}\quad
1<p<{\bf q}(n,\delta,k-1).
$$
$$
\hp{k}(M)=0 \quad\textrm{pour}\quad
1<p\leq{\bf q}(n,\delta,k).
$$}
\end{stheo}

\ali

Ce r{\'e}sultat, annonc{\'e} dans \cite{P'}, est un raffinement de celui de H. Donnelly et F. Xavier, \cite{Donnelly-Xavier}, concernant l'annulation de la
cohomologie $L^2$ r\'eduite. La condition d'annulation de la torsion est optimale. D'abord, pour l'espace hyperbolique ($\delta=-1$) en tout degr{\'e}, voir en \ref{hyperbreeltor}. Mais il y a d'autres exemples. Soient $n$ et $\mu$ des entiers tels que $2\leq\mu\leq n-1$ et $\delta\in]-1,0[$. Soit $G_{\mu,n,\delta}$ le produit semi-direct $G=\R\ltimes{\alpha}\R^{n-1}$ o{\`u} $\alpha$ est une matrice diagonale avec seulement deux valeurs propres
distinctes $1$ et $\sqrt{-\delta}<1$ de multiplicit{\'e}s $\mu-1$ et $n-\mu$. Alors le groupe de Lie $G_{\mu,n,\delta}$ poss{\`e}de une m{\'e}trique riemannienne invariante {\`a} gauche $\delta$-pinc{\'e}e.

\begin{stheo} {\bf B.~}
\label{B}
{\em Soient $n$ et $k=2,\ldots,n-1$ des entiers.
\begin{enumerate}
  \item Pour l'espace hyperbolique réel, $\tp{k}(\R H^n )\not=0$ si et seulement si $\d p=\frac{n-1}{k-1}$.
  \item Soient $\delta\in]-1,0[$ un r{\'e}el. Si $k=\mu$ et
$$
{\bf
q}(n,\delta,k-1)<p<1+{{1+(n-1-k)\sqrt{-\delta}}\over{k-2+\sqrt{-\delta}}},
$$
alors
$\tp{k}(G_{\mu,n,\delta})\not=0$, i.e. $\hp{k}(G_{\mu,n,\delta})$ n'est
pas s{\'e}par{\'e}.
\end{enumerate}

Par cons{\'e}quent, pour tout $\mu=2,\cdots,n-1$,
$G_{\mu,n,\delta}$ n'est pas quasiisom{\'e}trique {\`a} une
vari{\'e}t{\'e} $\delta'$-pinc{\'e}e avec $\delta'<\delta$.}
\end{stheo}

Ce résultat, qui élabore sur \cite{KS}, a été annoncé dans \cite{PCambridge}.

\subsection{Cas des espaces sym{\'e}triques de rang un}

Les espaces sym{\'e}triques de rang $1$ de type non compact {\`a} courbure non constante sont $-1/4$-pinc{\'e}s. Alors que ce sont de bons candidats pour tester l'optimalit{\'e} du th{\'e}or{\`e}me A, (la preuve ne comporte aucune perte quand on l'applique {\`a} ces espaces pour les valeurs ad{\'e}quates de $k$), le calcul r{\'e}v{\`e}le que leur cohomologie $L^p$ reste s{\'e}par{\'e}e au-del{\`a} des intervalles donn{\'e}s par le th{\'e}or{\`e}me A. Cela r{\'e}sulte de la non commutativit{\'e} de leur unipotent maximal, voir \cite{PP}.

\subsection{M{\'e}thode}

Une variété riemannienne $M$ à courbure sectionnelle négative ressemble à un produit. En effet, le flot $\phi_t$ de l'opposé du gradient d'une fonction de Busemann $b$ réalise un difféomorphisme de $M$ sur $H\times\R$, où $H=b^{-1}(0)$ est une horosphère. Par exemple, pour l'espace homogène $G_{2,4,-\inv4}$ (resp. le plan hyperbolique complexe $\C H^2$), $b(t)=-t$, $H=\R^3$ (resp. $H$ = groupe d'Heisenberg). Les orbites du flot $\phi_{t}$ sont des géodésiques asymptotes en $+\infty$, i.e. aboutissant en un même point à l'infini, en provenance de tous les autres points à l'infini (voir figure).

\begin{center}
\includegraphics[width=3in]{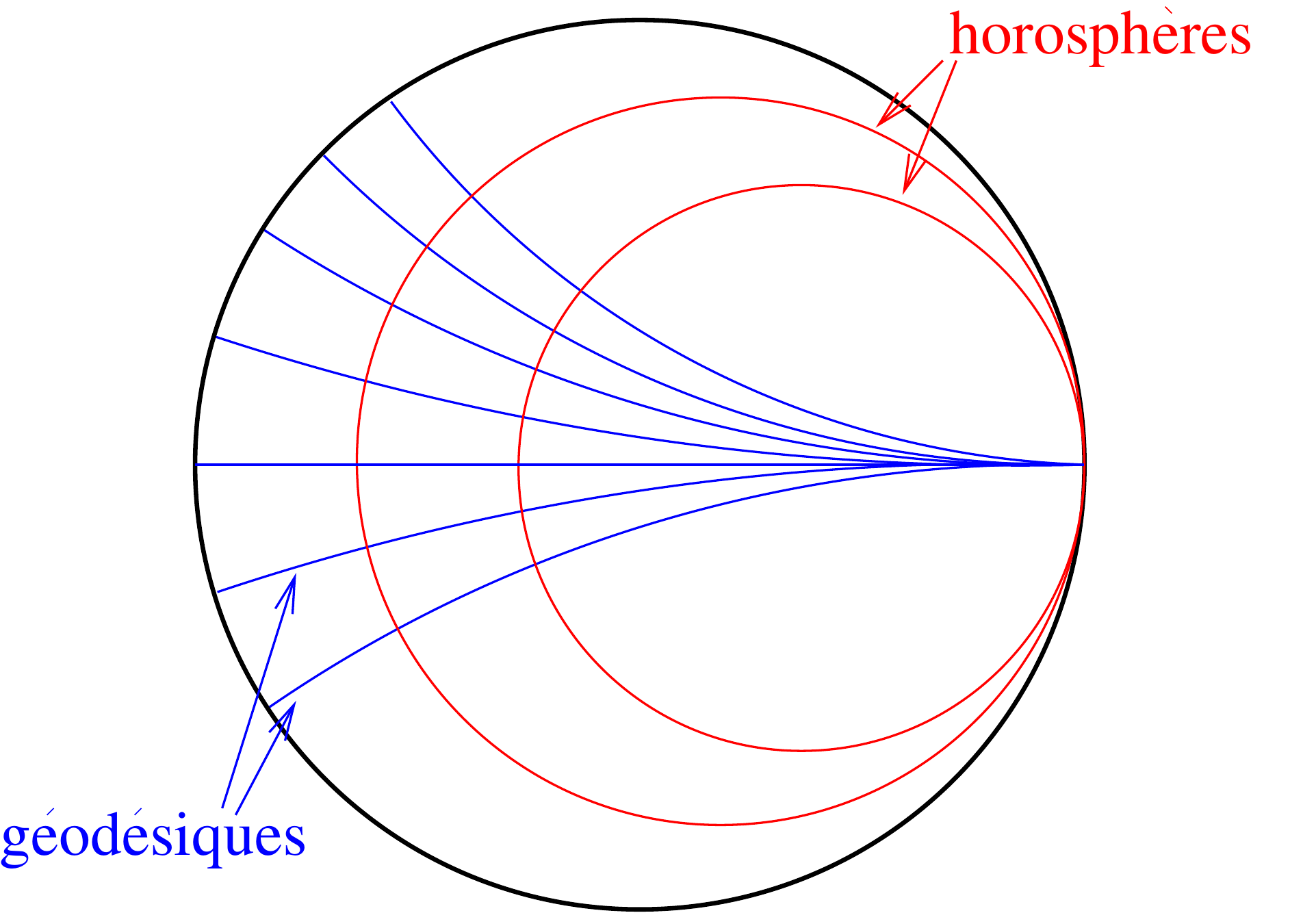}
\end{center}

On montre que si la courbure est suffisamment pincée (i.e., sous les hypothèses du théorème A), toute $k$-forme fermée $L^p$ $\omega$ possède une \emph{valeur au bord}
\begin{eqnarray*}
\omega_{\infty}=\lim_{t\to+\infty}\phi_{t}^{*}\omega.
\end{eqnarray*}
De plus, $\omega$ est la différentielle d'une forme $L^p$ si et seulement si $\omega_{\infty}=0$. Par conséquent, l'application valeur au bord induit une injection de $\hp{k}$ dans un espace séparé, donc $\hp{k}$ est séparé.

Inversement, pour les espaces homogènes $G_{k,n,\delta}$, on construit explicitement des classes de cohomologie non nulles, en utilisant la structure de produit semi-direct. Il faut se méfier de la formule de Künneth, qui n'est pas vraie en présence de torsion, même pour les produits directs. Après des préliminaires (dualité de Poincaré, annulation de la cohomologie $L^p$ réduite des groupes abéliens), on introduit et on construit des classes de torsion non nulles particulières, dites \emph{robustes}, qui restent non nulles après produit cartésien. La nature semi-directe du produit $G_{k,n,\delta}=\R^{n-1}\timez_{\alpha}\R$ exige la construction de classes robustes adaptées à la graduation de l'algèbre extérieure de $\R^{n-1}$ par les espaces propres de la dérivation $\alpha$. Puis on effectue le produit cartésien de ces classes avec des classes de cohomologie à support compact de $\R$. On obtient ainsi un intervalle ouvert de valeurs de $p$ pour lesquelles $\tp{k}(G_{k,n,\delta})\not=0$. Pour l'espace hyperbolique réel, il y a exactement une valeur de $p$ en chaque degré $>1$ pour laquelle $\tp{k}(\R H^n )\not=0$. On le montre en effectuant le produit cartésien d'une classe de torsion robuste de $\R$ avec une classe de cohomologie à support compact de $\R^{n-1}$.

\subsection{Remerciements}

Je tiens {\`a} remercier D. Rugina pour les nombreuses conversations que
nous avons eues autour de la cohomologie $L^{p}$, V. Goldshtein et M. Troyanov, pour leurs marques d'int{\'e}r{\^e}t et leur manuscrit \cite{GT} qui a {\'e}t{\'e} une source d'inspiration. 

\section{Annulation de la torsion}
\label{tor=0}

\subsection{Fonctions de Busemann}

Soit $M$ une variété riemannienne simplement connexe à courbure sectionnelle négative pincée. On se donne une fonction de Busemann $b$. C'est une fonction, obtenue comme limite de distances à des points, qui possède les propriétés suivantes.

\begin{enumerate}
  \item $b$ est lisse, son gradient est partout de norme 1.
  \item Les lignes de gradient de $b$ sont des géodésiques convergeant en $+\infty$ vers un même point du bord à l'infini de $M$.
  \item Les propriétés de contraction du flot $\phi_t$ de $-\nabla b$ sont contrôlées par la courbure sectionnelle. 
\end{enumerate}

\begin{exemple}
\label{buserh}
Cas de l'espace hyperbolique réel.
\end{exemple}
Dans ce cas, tout plan totalement géodésique contenant une ligne de gradient de $b$ est stable par $\phi_t$. Orthogonalement à ses orbites, $\phi_t$ est une homothétie de rapport $e^{-t}$ : $(\phi_t )^* (g-b^2 )=e^{-2t}(g-b^2 )$. Autrement dit, $\phi_t$ contracte de la même façon dans toutes les directions autour d'une orbite. $\phi_t$ multiplie les volumes par le facteur $e^{-(n-1)t}$, où $n=\dim M$. Si $\omega$ est une $k$-forme différentielle sur $M$, elle se décompose uniquement en $\omega=\beta+\gamma\wedge db$ de sorte que $\iota_\xi \beta=0$, $\iota_\xi \gamma=0$. Alors
\begin{eqnarray*}
|\phi_{t}^{*}\beta|(x)=e^{kt}|\beta|(\phi_t (x)),\quad |\phi_{t}^{*}\gamma|(x)=e^{(k-1)t}|\gamma|(\phi_t (x)).
\end{eqnarray*}
La formule de changement de variables donne
\begin{eqnarray*}
\int_{M}|\phi_{t}^{*}\beta|^p =e^{(kp-n+1)t}\int_{M}|\beta|^p ,\quad \int_{M}|\phi_{t}^{*}\gamma|^p =e^{((k-1)p-n+1)t}\int_{M}|\gamma|^p .
\end{eqnarray*}
Autrement dit, le flot $\phi_t$ contracte ou dilate exponentiellement la norme $L^p$ des $k$-formes différentielles, transversalement à ses orbites, suivant que $p$ est inférieur ou supérieur à $\frac{n-1}{k}$.

\begin{exemple}
\label{busech}
Cas de l'espace hyperbolique complexe. 
\end{exemple}
Dans ce cas, toute ligne de gradient est contenue dans une droite complexe, totalement géodésique, de courbure sectionnelle $-1$, stable par $\phi_t$. Tangentiellement à cette droite, et orthogonalement à l'orbite, $\phi_t$ est une homothétie de rapport $e^{-t}$. Tout plan contenant $\xi$ mais orthogonal à la droite complexe, s'exponentie en une surface totalement géodésique à courbure sectionnelle $-1$, stable par $\phi_t$, donc, dans ces directions, $\phi_t$ est une homothétie de rapport $e^{-t/2}$. Par conséquent, $\phi_t$ multiplie les volumes par $e^{-mt}$, où $m=\textrm{dim}_\C M=\frac{1}{2}\dim M$. Si $\omega$ est une $k$-forme différentielle sur $M$, elle se décompose uniquement en $\omega=\beta+\gamma\wedge db$ de sorte que $\iota_\xi \beta=0$, $\iota_\xi \gamma=0$, puis $\beta$ se décompose à son tour en $\beta=\epsilon+\eta\wedge Jdb$, où $J$ désigne la structure complexe, et $\iota_{J\xi}\epsilon=0$, $\iota_{J\xi}\eta=0$. Alors
\begin{eqnarray*}
|\phi_{t}^{*}\epsilon|(x)=e^{kt/2}|\epsilon|(\phi_t (x)),\quad |\phi_{t}^{*}\eta|(x)=e^{(k-1)t/2}|\gamma|(\phi_t (x)).
\end{eqnarray*}
Remarquer que $|\phi_{t}^{*}Jdb|=e^{-t}$. Il s'ensuit que
\begin{eqnarray*}
\int_{M}|\phi_{t}^{*}\epsilon|^p =e^{(kp-2m)t/2}\int_{M}|\epsilon|^p ,\quad \int_{M}|\phi_{t}^{*}\eta\wedge Jdb|^p =e^{((k+1)p-2m)t/2}\int_{M}|\eta\wedge Jdb|^p .
\end{eqnarray*}
Autrement dit, le flot $\phi_t$ contracte (resp. dilate) exponentiellement la norme $L^p$ de toutes les $k$-formes différentielles, transversalement à ses orbites, si $p<\frac{2m}{k+1}=\frac{n}{k}$ (resp. si $p>\frac{2m}{k}$). Lorsque $\frac{2m}{k+1}<p<\frac{2m}{k}$, la situation mérite plus d'attention.

\subsection{Champs de vecteurs $(k,p)$-contractants}

Les exemples ci-dessus suggèrent la définition suivante.

\begin{defi}
\label{kpcontractant}
Soit $M$ une variété riemannienne. Soit $\xi$ un champ de vecteurs complet sur $M$, soit $\phi_t$ son flot. Soit $p>1$ un réel, soit $k$ un entier inférieur à la dimension de $M$. On dit que $\xi$ est \emph{$(k,p)$-contractant} si $\phi_t$ diminue exponentiellement la norme $L^p$ des $k$-formes transversalement à $\xi$. Plus précisément, on note $Jac(\phi_t )$ le jacobien de $\phi_t$, et on demande qu'il existe des constantes $C$ et $\eta>0$ telles que, pour tout $x\in M$ et toute $k$-forme $\beta\in\Lambda^k T^* M$ telle que $\iota_\xi \beta=0$,
\begin{eqnarray*}
|\phi_{t}^{*}\beta|(x)\,Jac_{x}(\phi_t )^{1/p}\leq C\,e^{-\eta t}|\beta|(\phi_t (x))
\end{eqnarray*}
pour tout $t\geq 0$.

On dit que $\xi$ est \emph{$(k,p)$-dilatant} si $-\xi$ est $(k,p)$-contractant.
\end{defi}

\begin{exemple}
Cas des produits semi-directs $G=H\timez_{\alpha}\R$.
\end{exemple}
Ici, $H$ est un groupe de Lie, $\alpha$ une dérivation de l'algèbre de Lie de $H$ qui engendre un groupe à un paramètre $e^{t\alpha}$ d'automorphismes de $H$, et $G=H\times\R$ muni de la multiplication
\begin{eqnarray*}
(h,t)(h',t')=(h\,e^{t\alpha}(h'),t+t').
\end{eqnarray*}
On utilise le champ de vecteurs invariant à gauche $\xi=\ddt$ qui engendre l'action à droite du facteur $\R$. Alors les formes différentielles annulées par $\iota_{\xi}$ s'identifient aux formes différentielles sur $H$ dépendant de $t$. Notons $sp(\alpha)$ l'ensemble des valeurs propres de $\alpha$ répétées autant de fois que leurs multiplicités. Le flot $\phi_t$ agit sur les $k$-formes transverses avec pour valeurs propres les nombres $e^{-t\lambda}$, où $\lambda$ décrit les sommes de $k$ éléments de $sp(\alpha)$. Par conséquent, $\xi$ est $(k,p)$-contractant si et seulement si les parties réelles de toutes ces sommes sont strictement supérieures à $\frac{\tr(\alpha)}{p}$.

\begin{prop}
\label{pince}
Soit $M$ une vari{\'e}t{\'e} riemannienne compl{\`e}te de dimension $n$,
simplement connexe, dont la courbure sectionnelle $K$ satisfait $-1\leq K\leq\delta<0$. Soit $\xi$ un champ de vecteurs de Busemann. Si $k=0,\cdots,n-1$ et si $p>1$ satisfait
$$
p<\mathbf{q}(n,\delta,k)=1+{{n-k-1}\over{k}}\sqrt{-\delta},\quad\hbox{ (resp. }
p>1+{{n-k-1}\over{k\sqrt{-\delta}}}),
$$
alors le champ $\xi$ est $(k,p)$-contractant (resp. $(k,p)$-dilatant).
\end{prop}

\preuve
Notons $\phi_t$ le flot de $\xi$. Ses trajectoires sont des g{\'e}od{\'e}siques
parcourues {\`a} vitesse $1$. Soit $x\in M$. La quantit{\'e} {\`a} majorer est
$$
n(t,x)=p\log\n{(\L^{k}d\phi_{t})_{|\L^{k}\xi^{\perp}}}-\log det(d\phi_{t}).
$$
Elle satisfait, pour tous $s$ et $t$, $n(t+s,x)=n(s,\phi_{t}(x))+n(t,x)$.

Notons $\tau_{t}$ le transport parall{\`e}le de $\phi_{t}(x)$ {\`a} $x$ le long de la g{\'e}od{\'e}sique $s\mapsto \phi_{s}(x)$. Alors $\tau_{t}d\phi_{t}$
pr{\'e}serve l'hyperplan orthogonal {\`a} $\xi(x)$. Notons $J(t)$ sa matrice
dans une base orthonorm{\'e}e de $\xi(x)^{\perp}$, de sorte que
$\d n(t,x)=p\log\n{\L^{k}J(t)}-\log det(J(t))$.

Comme $\xi$ est un gradient, la matrice $U(t)=J(t)^{-1}J'(t)$, seconde forme
fondamentale des hypersurfaces de niveau, est sy\-m{\'e}\-tri\-que. Com\-me les colonnes de $J$  sont des champs de Jacobi, la matrice $U(t)$ satisfait l'{\'e}quation de Riccati
$$
U'+U^{2}+R=0
$$
o{\`u} $R$ est la matrice de l'op{\'e}rateur de courbure $v\mapsto
R(v,\xi)\xi$. Classiquement (voir par exemple \cite{courbure}, \cite{CE}), on en tire une estimation des valeurs propres $\lambda_{1},\ldots,\lambda_{n-1}$ de $U$,
$$
\sqrt{-\delta}\leq \lambda_{1}\leq\cdots\leq\lambda_{n-1} \leq 1.
$$
Comme $J(0)=I$ est l'identit{\'e}, $J(t)=I+tU(0)+o(t)$ donc $\n{\L^{k}J(t)}\leq
1+|t|\n{{\cal D}^{k}U(0)}+o(t)$ o{\`u} ${\cal D}^{k}U$ d{\'e}signe l'extension de $U$ comme d{\'e}\-ri\-va\-ti\-on de l'alg{\`e}bre ext{\'e}rieure. On peut donc majorer la d{\'e}riv{\'e}e {\`a} droite
\begin{eqnarray*}
n'(0+)={{\partial n}\over{\partial t}}(0,x)
&\leq& p\n{{\cal D}^{k}U(0)}-{\rm tr U(0)}\\
&\leq& p(\sum_{i=n-k}^{n-1} \l_{i})-\sum_{i=1}^{n-1} \l_{i}\\
&=& (p-1)(\sum_{i=n-k}^{n-1} \l_{i})-\sum_{i=1}^{n-k-1} \l_{i}\\
&\leq& k(p-1)-(n-k-1)\sqrt{-\delta}.
\end{eqnarray*}
En d{\'e}rivant l'{\'e}quation $n(t+s,x)=n(s,\phi_{t}(x))+n(t,x)$, on trouve que $n'(t+,x)=n'(0+,\phi_{t}(x))\leq k(p-1)-(n-k-1)\sqrt{-\delta}$ pour tout
$t$. En int{\'e}grant, il vient pour tout $t\in\R$,
$$
\n{(\L^{k}d\phi_{t})_{|\L^{k}\xi^{\perp}}}^{p}\leq e^{-\eta t}Jac(\phi_{t}),
$$
avec $\eta=(n-k-1)\sqrt{-\delta}-k(p-1)$. Si $\eta>0$, i.e. si la courbure
est suffisamment pinc{\'e}e, on conclut que $\xi$ est $(k,p)$-contractant.

Si on remplace $\xi$ par $-\xi$, les valeurs propres $\l_i$ de la seconde forme fondamentale sont remplac{\'e}es par $\mu_{i}=-\l_{n-i}$ qui satisfont
$$
-1\leq \mu_{1}\leq\cdots\leq\mu_{n-1} \leq -\sqrt{-\delta}.
$$
La nouvelle fonction ${\tilde n}(t,x)=n(-t,x)$ satisfait
\begin{eqnarray*}
{\tilde n}'(0+)
&\leq& p(\sum_{i=n-k}^{n-1} \mu_{i})-\sum_{i=1}^{n-1} \mu_{i}\\
&=& (p-1)(\sum_{i=n-k}^{n-1} \mu_{i})-\sum_{i=1}^{n-k-1} \mu_{i}\\
&\leq& k(p-1)(-\sqrt{-\delta})+(n-k-1).
\end{eqnarray*}
Il vient
$$
\n{(\L^{k}d\phi_{t})_{|\L^{k}\xi^{\perp}}}^{p}\leq e^{\eta' t}Jac(\phi_{t})
$$
avec $\eta'=k(p-1)\sqrt{-\delta}-n+k+1$. Si $\eta'>0$, on conclut que $\xi$ est $(k,p)$-dilatant.\qed

\begin{rem}
\label{limite}
Cas limite.
\end{rem}
Si $p=\mathbf{q}(n,\delta,k)$, le flot $\phi_t$ diminue au sens large la norme $L^p$ des $k$-formes transverses, au sens où 
$$
\n{(\L^{k}d\phi_{t})_{|\L^{k}\xi^{\perp}}}^{p}\leq Jac(\phi_{t}).
$$

\begin{rem}
\label{egalite}
Cas d'égalité.
\end{rem}
Dans l'argument ci-dessus, les in{\'e}galit{\'e}s sont optimales dans le cas où les valeurs propres ne prennent que deux valeurs. Il est facile, à l'aide de \cite{Heintze}, de faire la liste des espaces homog{\`e}nes {\`a} courbure sectionnelle strictement n{\'e}gative pour lesquels les valeurs propres prennent exactement deux valeurs {\'e}gales aux bornes de la courbure sectionnelle. En voici deux familles particulières.

\begin{exemple}
Les espaces sym{\'e}triques de rang un. 
\end{exemple} 
La courbure sectionnelle varie entre $-1$ et $-1/4$. Les valeurs propres sont $1/2$ (avec multiplicit{\'e} $2m-2$ pour l'espace hyperbolique complexe ${{\bf C}}H^{m},~m\geq 2$, $4m-4$ pour l'espace hyperbolique quaternionien ${{\bf H}}H^{m},~m\geq 2$, $8$ pour le plan hyperbolique des octaves de Cayley ${{\bf Ca}}H^{2}$), et $1$ avec multiplicit{\'e} compl{\'e}mentaire, soit respectivement 1, 3 et 7.

\begin{exemple}
Une famille d'espaces homogènes. 
\end{exemple} 
Soient $1\leq\mu<n$ des entiers et $\delta\in]-1,0[$. Soit $G_{\mu,n,\delta}$ le produit semi-direct $G=\R^{n-1}\timez_{\alpha}\R$ o{\`u} $\alpha$ est une matrice diagonale avec seulement deux valeurs propres distinctes $1$ et $\sqrt{-\delta}<1$ de multiplicit{\'e}s $\mu-1$ et $n-\mu$. La m{\'e}trique invariante $dt^{2}+e^{2t}dx^{2}+e^{2t\sqrt{-\delta}}dy^{2}$ (o{\`u} $x$ regroupe les $\mu-1$ premi{\`e}res coordonn{\'e}es de $\R^{n-1}$ et $y$ les $n-\mu$ suivantes) a une courbure sectionnelle comprise entre $-1$ et $-\delta$.

\subsection{Valeur au bord}

\begin{prop}
\label{valeur}
Soit $M$ une variété riemannienne, soit $\xi$ un champ de vecteurs complet sur $M$, de flot $\phi_t$. 

\begin{enumerate}
  \item On suppose que $\xi$ est $(k-1,p)$-contractant et que sa norme est bornée. Alors toute $k$-forme fermée $L^p$ $\omega$ possède une \emph{valeur au bord}
\begin{eqnarray*}
\omega_{\infty}=\lim_{t\to+\infty}\phi_{t}^{*}\omega,
\end{eqnarray*}
et $\omega-\omega_{\infty}$ est la différentielle d'une forme $L^p$.
  \item Si $\xi$ est $(k-1,p)$- et $(k-2,p)$-contractant, et si $\omega=d\beta$ où $\beta\in L^p$, alors $\omega_{\infty}=0$.
  \item Si $\xi$ est $(k-1,p)$- et $(k,p)$-contractant, alors $\omega_{\infty}=0$.
\end{enumerate} 
Par conséquent, 
\begin{enumerate}
  \item Si $\xi$ est $(k-1,p)$- et $(k-2,p)$-contractant, $\tp{k}(M)=0$.
  \item Si $\xi$ est $(k-1,p)$- et $(k,p)$-contractant, $\hp{k}(M)=0$.
\end{enumerate} 
\end{prop}

\preuve
D'après la formule de Cartan, la dérivée de Lie $\mathcal{L}_{\xi}\omega=\ddt \phi_{t}^{*}\omega_{|t=0}$ est égale à
\begin{eqnarray*}
\mathcal{L}_{\xi}\omega=d(\iota_{\xi}\omega)+\iota_{\xi}(d\omega).
\end{eqnarray*}
Supposons que $d\omega=0$. En intégrant l'identité $\ddt \phi_{t}^{*}\omega=\phi_{t}^{*}\mathcal{L}_{\xi}\omega$, il vient
\begin{eqnarray*}
\phi_{t}^{*}\omega-\omega&=&\int_{0}^{t}\phi_{s}^{*}\mathcal{L}_{\xi}\omega\,ds\\
&=&d(\int_{0}^{t}\phi_{s}^{*}\iota_{\xi}\omega\,ds).
\end{eqnarray*}
Si $\xi$ est borné, $\n{\iota_{\xi}\omega}\lp\leq\n{\xi}_{L^{\infty}}\n{\omega}\lp$. Si de plus $\xi$ est $(k-1,p)$-contractant, il existe $C$ et $\eta>0$ tels que $\n{\phi_{s}^{*}\iota_{\xi}\omega}\lp \leq C\,e^{\eta s}\n{\omega}\lp$. Par conséquent, l'intégrale
\begin{eqnarray*}
B\omega=\int_{0}^{+\infty}\phi_{s}^{*}\iota_{\xi}\omega\,ds
\end{eqnarray*}
converge dans $L^p$. On note
\begin{eqnarray*}
\omega_{\infty}&=&\omega+dB\omega\\
&=&\lim_{t\to +\infty}\phi_{t}^{*}\omega.
\end{eqnarray*}
Il s'agit d'une limite au sens des distributions. Si la limite est nulle, alors $\omega=d(-B\omega)$ au sens des distributions. Cela entraîne que $-B\omega\in\Omega^{k-1,p}(M)$, et que sa différentielle est $\omega$, donc que la classe de cohomologie $L^p$ de $\omega$ est nulle.

Si $\xi$ est de plus $(k-2,p)$-contractant, on peut aussi définir un opérateur $B$ borné sur les formes $L^p$ de degré $k-1$. Soit $\alpha$ une $k-1$-forme $L^p$ telle que $d\beta=\omega$. Il vient
\begin{eqnarray*}
\phi_{t}^{*}\beta-\beta&=&\int_{0}^{t}\phi_{s}^{*}\mathcal{L}_{\xi}\omega\,ds\\
&=&d(\int_{0}^{t}\phi_{s}^{*}\iota_{\xi}\beta\,ds)+
\int_{0}^{t}\phi_{s}^{*}\iota_{\xi}(d\beta)\,ds,
\end{eqnarray*}
qui tend vers $dB\beta+B\omega$ quand $t$ tend vers $+\infty$. Mais comme $\xi$ est $(k,p)$-contractant, $\phi_{t}^{*}\beta$ tend vers 0. On trouve que $\beta=-dB\beta-B\omega$, d'où 
\begin{eqnarray*}
\omega=d\beta=-dB\omega=\omega-\omega_{\infty},
\end{eqnarray*}
d'où $\omega_{\infty}=0$.

Cela prouve que $dL^p$ est exactement le noyau de l'application valeur au bord, de $\Omega^{k,p}(M)\cap\ker d$ dans l'espace vectoriel topologique des formes différentielles sur $M$ à coefficients distributions. Par conséquent, il est fermé, donc $\tp{k}(M)=0$.

Supposons que $\xi$ est $(k-1,p)$- et $(k,p)$-contractant. Soit $\omega$ une $k$-forme fermée $L^p$. On écrit $\omega=\beta+db\wedge\gamma$ où $\iota_{\xi}\beta=0$, $\iota_{\xi}\gamma=0$. Alors $\phi_{t}^{*}\beta$ et $\phi_{t}^{*}\gamma$ tendent vers 0 dans $L^p$, donc $\phi_{t}^{*}\omega$ tend vers 0 dans $L^p$, donc $\omega_{\infty}=0$, d'où $\omega\in dL^p$. Cela prouve que $\hp{k}(M)=0$.\qed

\begin{rem}
Plus généralement, sous des hypothèses adéquates, l'opérateur $B$ définit une homotopie du complexe $\Omega^{*,p}(M)$ sur un complexe de formes différentielles invariantes par le flot $\phi_t$. Ce point de vue est développé dans \cite{GKS}, \cite{PP}. 
\end{rem}

\subsection{Preuve du théorème A}

Soit $M$ une vari{\'e}t{\'e} riemannienne compl{\`e}te, simplement connexe, à courbure négative $\delta$-pincée. Soit $k<n=\dim M$. Notons $\d{\bf q}(n,\delta,k)=1+{{n-k-1}\over{k}}\sqrt{-\delta}$. Remarquer que ${\bf q}(n,\delta,k)$ est une fonction décroissante de $k$.

D'après la proposition \ref{pince}, si $p<{\bf q}(n,\delta,k-1)$, les champs de vecteurs de Busemann $\xi$ sont $(k-1,p)$ et $(k-2,p)$-contractants. La proposition \ref{valeur} s'applique, et $\tp{k}(M)=0$. De même, si $p<{\bf q}(n,\delta,k)$, $\xi$ est $(k-1,p)$ et $(k,p)$-contractant, donc $\hp{k}(M)=0$.

Il reste à traiter le cas limite $p={\bf q}(n,\delta,k)$. Dans ce cas, d'après la remarque \ref{limite},
$$
\label{limi}
\n{(\L^{k}d\phi_{t})_{|\L^{k}\xi^{\perp}}}^{p}\leq Jac(\phi_{t}).
$$
Soit $K$ un compact de $M$. Il existe une constante $c=c(K)$ telle que les images $\phi_{cj} (K)$ pour $j\in\Z$ soient deux à deux disjointes. Alors la suite $\n{\omega}_{L^p (\phi_{cj} (K))}$ est dans $\ell^p (\Z)$, donc tend vers 0. L'inégalité \ref{limi} entraîne que si $\omega$ est une $k$-forme sur $M$ annulée par $\iota_{\xi}$,
\begin{eqnarray*}
\n{\phi_{cj}^{*}\omega}_{L^p (K)}\leq\n{\omega}_{L^p (\phi_{cj} (K))}
\end{eqnarray*}
qui tend vers 0. Cela montre que la limite au sens des distributions $\omega_{\infty}$ est nulle sur tout compact, donc est nulle. On conclut que $\hp{k}(M)=0$ aussi dans ce cas.\qed

\begin{rem}
Cas des produits semi-directs $G=H\timez_{\alpha}\R$.
\end{rem}
Dans le cas des groupes $G_{\mu,n,\delta}$, le théorème A s'applique, et la torsion $L^p$ s'annule en degré $\mu$ pour tout $p<1+{\bf q}(n,\delta,\mu-1)$.

Soit $G=H\timez_{\alpha}\R$ un produit semi direct plus général. Notons $\lambda_1 \leq \cdots\leq \lambda_{n-1}$ les parties réelles des valeurs propres de $\alpha$ répétées autant de fois que leur multiplicité. On suppose que $\lambda_1 >0$.  On utilise le champ de vecteurs invariant à gauche $\xi$ qui engendre l'action à droite du facteur $\R$. Le champ $-\xi$ est $(k-1,p)$-contractant et $(k-2,p)$-contractant tant que $p$ reste strictement inférieur à $\d\frac{\tr\alpha}{\lambda_1 +\cdots+\lambda_{k-1}}$. La proposition \ref{valeur} s'applique, et on conclut que la torsion $L^p$ s'annule. 

\section{Exemples où la torsion est non nulle}

Ce seront des groupes de Lie, produits semi-directs de groupes abéliens par $\R$. On construit des formes différentielles fermées explicites en utilisant la structure produit. Elles sont nulles en cohomologie réduite, parce que la cohomologie réduite d'un groupe nilpotent est nulle. Pour montrer qu'elles sont non nulles en torsion $L^p$, on utilise la dualité de Poincaré.

\subsection{Dualit{\'e} de Poincar{\'e}}
\label{poin}

Le lemme suivant est essentiellement d{\^u} {\`a} V. Goldshtein et M. Troyanov, \cite{GT}.

\begin{lemme}
\label{stokes}
Soit $M$ une vari{\'e}t{\'e} riemannienne orient{\'e}e compl{\`e}te de dimension $n$.
Etant donn{\'e} $p>1$, on note
$p'$ l'exposant conjugu{\'e}, i.e. tel que $\d\inv{p}+\inv{p'}=1$. Soit
$\omega$ une $k$-forme diff{\'e}rentielle ferm{\'e}e et $L^p$ sur $M$. Alors
\begin{itemize}
\item $\omega$ est non nulle en cohomologie $L^p$ r{\'e}duite si et seulement si
il existe une $n-k$-forme ferm{\'e}e $\psi\in L^{p'}$ telle que
$\int_{M}\omega\wedge\psi\not=0$.

\item $\omega$ est non nulle en cohomologie $L^p$ si et seulement si
il existe une suite $\psi_j$ de $n-k$-formes diff{\'e}rentielles $L^{p'}$ telles
que $\int_{M}\omega\wedge\psi_{j}\geq 1$ et $\n{d\psi_{j}}_{L^{p'}}$ tend vers
$0$.
\end{itemize}
\end{lemme}

\preuve
Comme $M$ est compl{\`e}te, pour toute $n-1$-forme $L^1$ dont la
diff{\'e}ren\-ti\-el\-le est $L^1$, on a
$$
\int_{M}d\omega=0.
$$
Par cons{\'e}quent, si $\omega\in \op{k}(M)$ et $\psi\in\Omega^{n-1-k,p'}(M)$,
$$
\int_{M}\omega\wedge d\psi=(-1)^{k+1}\int_{M} d\omega\wedge\psi.
$$

Si $\omega\in\op{k}(M)$ est nulle en cohomologie $L^p$ r{\'e}duite, alors il
existe une suite $\beta_{j}\in\op{k-1}(M)$ telle que $d\beta_{j}$
converge vers $\omega$ dans $L^p$. Si $\psi\in\Omega^{n-1-k,p'}(M)$, il vient
$$
\int_{M} \omega\wedge\psi=\lim_{j} \int_{M} d\beta_{j}\wedge \psi
=\lim_{j} \int_{M}\beta_{j}\wedge d\psi.
$$
Par cons{\'e}quent, pour toute forme ferm{\'e}e $\psi\in\Omega^{n-k,p'}(M)$,
$\int_{M}\omega\wedge\psi=0$. 

Inversement, si $\omega\in\op{k}(M)$ n'est pas nulle en cohomologie r{\'e}duite,
alors, d'apr{\`e}s Hahn-Banach, il existe une forme lin{\'e}aire continue $L$ sur
$L^{p}\Omega^{k}(M)$ qui s'annule sur l'adh{\'e}rence de l'image $d\op{k-1}(M)$
mais pas sur $\omega$. Par dualit{\'e} de $L^p$ et
$L^{p'}$, il existe une forme $\psi\in L^{p'}\Omega^{n-k}(M)$ telle
que pour tout $\gamma\in L^{p}\Omega^{k}(M)$,
$L(\gamma)=\int_{M}\gamma\wedge\psi$. Si $\beta$ est lisse et {\`a} support
compact, on a $0=L(d\beta)=\int_{M} d\beta\wedge\psi$, i.e. $d\psi=0$ au sens
des distributions. On conclut que $\psi\in \Omega^{n-1-k,p'}(M)$ est ferm{\'e}e
et satisfait $\int_{M}\omega\wedge\psi\not=0$.

Si $\omega\in\op{k}(M)$ est nulle en cohomologie $L^p$, i.e. $\omega=d\beta$
o{\`u} $\beta\in\op{k-1}(M)$, alors pour tout $\psi\in \Omega^{n-1-k,p'}(M)$, 
$$
\int_{M}\omega\wedge\psi=\int_{M}\beta\wedge d\psi\leq\n{\beta}_{L^{p}}
\n{d\psi}_{L^{p'}},
$$
donc si $\d \n{d\psi_{j}}_{L^{p'}}$ tend vers $0$, il en est de m{\^e}me de $\d
\int_{M}\omega\wedge\psi_{j}$.

Inversement, soit $\omega\in\op{k}(M)$ une forme ferm{\'e}e. Si $\omega$ n'est
pas nulle en cohomologie r{\'e}duite, il existe une $n-k$-forme ferm{\'e}e $\psi\in L^{p'}$ telle que
$\int_{M}\omega\wedge\psi\not=0$. La suite stationnaire $\psi_{j}=\psi$ pour
tout $j$ convient. Supposons d{\'e}sormais que $\omega$ est nulle en cohomologie
r{\'e}duite. On d{\'e}finit une forme lin{\'e}aire $L$ sur $d\Omega^{n-k,p'}(M)$
comme suit. Etant donn{\'e} $\gamma\in d\Omega^{n-k,p'}(M)$, on choisit
$\psi\in \Omega^{n-k,p'}(M)$ tel que $d\psi=\gamma$ et on pose
$L(\gamma)=\int_{M}\omega\wedge\psi$. Comme l'int{\'e}grale de $\omega$ contre
une forme ferm{\'e}e est toujours nulle, le r{\'e}sultat ne d{\'e}pend pas du choix de
$\psi$. Supposons qu'il n'existe pas de suite
$\psi_{j}\in\Omega^{n-1-k,p'}(M)$ telle que $\int_{M}\omega\wedge\psi_{j}\geq
1$ et $\n{\psi_{j}}_{L^{p'}}$ tend vers $0$. 
Alors la forme lin{\'e}aire $L$ est continue pour la norme $L^{p'}$. Par
Hahn-Banach, $L$ se prolonge en une forme lin{\'e}aire continue sur
$L^{p'}\Omega^{n-k+1}(M)$. Par dualit{\'e} entre $L^p$ et $L^{p'}$, il existe une
$k-1$-forme $\beta\in L^p$ telle que pour tout $\gamma\in
L^{p'}\Omega^{n-k+1}(M)$, $L(\gamma)=(-1)^{k}\int_{M}\beta\wedge\gamma$. Si
$\psi$ est lisse {\`a} support compact, 
$$ 
\int_{M}\beta\wedge d\psi=(-1)^{k}\int_{M}\omega\wedge\psi
$$
donc $d\beta=\omega$ au sens des distributions. Par cons{\'e}quent,
$\beta\in\op{k-1}(M)$ et $\omega$ est nulle en cohomologie $L^p$.\qed

\begin{cor}
\label{poincare}
Soit $M$ une vari{\'e}t{\'e} riemannienne compl{\`e}te de dimension $n$. Soit $p>1$ et
$p'=p/(p-1)$. Alors 
$$
\rp{k}(M)\Leftrightarrow R^{n-k,p'}(M)=0,\quad \tp{k}(M)=0\Leftrightarrow
T^{n-k+1,p'}(M)=0.
$$
\end{cor}

\preuve
L'{\'e}nonc{\'e} sur la cohomologie r{\'e}duite r{\'e}sulte imm{\'e}diatement du lemme
\ref{stokes}.

Supposons que $T^{n-k+1,p'}(M)=0$. Montrons qu'il existe une constante $C$ telle que pour tout $\psi\in\Omega^{n-k,p'}(M)$, il existe $\gamma\in\Omega^{n-k,p'}(M)$ telle que $d\gamma=d\psi$ et $\n{\gamma}_{L^{p'}}\leq C\,\n{d\psi}_{L^{p'}}$. Par hypothèse, $d\Omega^{n-k,p'}(M)$ est ferm{\'e} dans $\Omega^{n-k+1,p'}(M)$. L'op{\'e}rateur $d$ induit $\overline{d}:\Omega^{n-k,p'}(M)/\ker d\to d\Omega^{n-k,p'}(M)$. C'est une bijection continue entre espaces de Banach, donc un isomorphisme. Notons $\overline{d}^{-1}$ son inverse, notons $C$ la norme de cet opérateur. Etant donn{\'e}e une $n-k$-forme $\omega\in\Omega^{n-k,p'}(M)$, soit $\phi\in\Omega^{n-k,p'}(M)$ un repr{\'e}sentant de la classe $\overline{d}^{-1}d\omega\in \Omega^{n-k,p'}(M)/\ker d$ de norme presque minimum. Elle satisfait (presque)
$$
d\omega=d\phi\et \n{\phi}_{L^{p'}} \leq C\,\n{d\phi}_{L^{p'}} .
$$

Soit
$\omega\in\op{k}(M)$ une forme ferm{\'e}e, nulle en cohomologie r{\'e}duite. Alors
$$
\int_{M}\omega\wedge\psi=\int_{M}\omega\wedge\gamma
$$
est contr{\^o}l{\'e}e par $\n{d\psi}_{L^{p'}}$. Par cons{\'e}quent, il n'existe pas
de suite $\psi_j$ de $n-k$-formes diff{\'e}rentielles $L^{p'}$ telles que
$\int_{M}\omega\wedge\psi_{j}\geq 1$ et $\n{d\psi_{j}}_{L^{p'}}$ tende vers $0$.
On conclut que $\omega$ est nulle en cohomologie $L^{p'}$. \qed

\remarque
Plus g{\'e}n{\'e}ralement, $\rp{k}(M)$ est isomorphe au dual de $R^{n-k,p'}(M)$. On aimerait dire que $\tp{k}(M)=Ext(T^{n-k+1,p'}(M),\R)$ dans une cat{\'e}gorie ad{\'e}quate.

\subsection{Cohomologie réduite des groupes abéliens}

On va construire des classes de cohomologie $L^p$ non nulles. Pour montrer qu'elles appartiennent à la torsion, nous auront besoin, au cours d'un raisonnement, de savoir que la cohomologie réduite de $\R^{n-1}$ est nulle.

La remarque suivante appara\^{\i}t entre autres dans \cite{Gromovasympt}. Elle s'applique notamment aux groupes de Lie nilpotents simplement connexes.

\begin{prop}
\label{rpnul}
Soit $G$ un groupe de Lie simplement connexe dont l'alg\`ebre de Lie a un centre non trivial. Alors la cohomologie $L^p$ r{\'e}duite de $G$ est nulle en tous degr{\'e}s.
\end{prop}

\preuve
Un vecteur non nul du centre donne un champ de vecteurs de Killing $\xi$ de longueur constante. La formule
$$
(\phi_{t})^{*}\omega-\omega=d\int_{0}^{t}(\phi_{s})^{*}\iota_{\xi}\omega\,ds
+\int_{0}^{t}(\phi_{s})^{*}\iota_{\xi}d\omega\,ds
$$
montre que le flot de $\xi$ agit trivialement sur la cohomologie $L^{p}$. Soit $\omega$ une forme ferm{\'e}e $L^{p}$ non nulle en cohomologie r\'eduite. Il existe donc une forme ferm{\'e}e $L^{p'}$ $\psi$ telle que $\d\int_{G}\omega\wedge\psi\not=0$. Comme le flot $\phi_{t}$ est l'identit\'e en cohomologie,
$$
\int_{G}(\phi_{t})^{*}\omega\wedge\psi=\int_{G}\omega\wedge\psi
$$
pour tout $t$. 

On utilise maintenant le fait que l'action de $\R$ sur $G$ par le flot de $\xi$ est propre. Soit $K$ un compact tel que la norme $L^{p}$ (resp. $L^{p'}$) de $\omega$ (resp. $\psi$) dans $G\setminus K$ soit petite. Soit $t$ tel que $\phi_{t}(K)$ soit disjoint de $K$. Alors 
$$
\int_{G}(\phi_{t})^{*}\omega\wedge\psi
\leq\n{\omega}_{L^{p}(G\setminus K)}\n{\psi}_{L^{p'}(G)}+
\n{\psi}_{L^{p'}(G\setminus K)}\n{\omega}_{L^{p}(G)}
$$
est petit, contradiction.\qed 

\subsection{Torsion des produits directs}
\label{prod}

L'objectif est d'étudier la torsion $L^p$ des groupes de Lie $G_{\mu,n,\delta}$. Il s'agit de produits semi-directs. Une première étape consiste à comprendre les produits directs, ou plus généralement, les produits riemanniens.

\ali

Soient $M_{1}$ et $M_{2}$ deux vari{\'e}t{\'e}s riemanniennes compl{\`e}tes. On note $\pi_i : M_1 \times M_2 \to M_i$ les projections. Lorsqu'elle est vraie, la formule de Künneth énonce que le produit cartésien des formes différentielles,
\begin{eqnarray*}
(\alpha_1 ,\alpha_2 )\mapsto \alpha_1 \times \alpha_2 =\pi_{1}^{*}\alpha_1 \wedge \pi_{2}^{*}\alpha_2 ,\quad \Omega^{*,p}(M_1 )\otimes \Omega^{*,p}(M_2 )\to\Omega^{*,p}(M_1 \times M_2 ),
\end{eqnarray*}
induit un isomorphisme en cohomologie $L^p$. Si la torsion $\tp{*}(M_{1})$ est identiquement nulle, c'est vrai, voir \cite{GKS}. Mais si $\tp{*}(M_{1})$ et $\tp{*}(M_{2})$ ne sont pas nulles, il en va autrement.

\begin{exemple}
Pour $p=2$, il existe des classes non nulles $\alpha$, $\beta\in T^{1,2}(\R)$ telles que $\alpha\times\beta=0$ dans $H^{2,2}(\R^2 )$.
\end{exemple}

\preuve
\def\g{\hat{\gamma}}
\def\a{\hat{a}}
\def\b{\hat{b}}
Soient $\alpha=a(x)\,dx$ une 1-forme $L^2$ et $f$ une fonction $L^2$ sur $\R$. L'équation $df=\alpha$, en Fourier, s'écrit
\begin{eqnarray*}
i\xi\hat{f}(\xi)=\a(\xi).
\end{eqnarray*}
Par conséquent, la classe de cohomologie $L^2$ de $\alpha$ est nulle si et seulement si $\d\xi\mapsto\frac{\a(\xi)}{i\xi}$ est $L^2$.

Etant données des 1-formes $L^2$ $\alpha=a(x)\,dx$, $\beta=b(y)\,dy$ sur $\R$ et 1-forme $L^2$ $\gamma=\gamma_x dx+\gamma_y dy$ sur $\R^2$, l'équation $d\gamma=\alpha\times\beta$ se traduit par
\begin{eqnarray*}
 i\xi\g_x (\xi,\eta)-i\eta\g_y (\xi,\eta)= \a(\xi)\b(\eta). 
\end{eqnarray*}
On la résoud en prenant
\begin{eqnarray*}
\g_x (\xi,\eta)=-\frac{i\xi}{\xi^2 +\eta^2}\a(\xi)\b(\eta),\quad \g_y (\xi,\eta)=\frac{i\eta}{\xi^2 +\eta^2}\a(\xi)\b(\eta).
\end{eqnarray*}
On choisit pour $\a=\b$ une fonction paire, lisse, à support compact qui, au voisinage de 0, coïncide avec $|\log(1/|\xi|)|^{-1/2}$. Alors $\a(\xi)/\xi$ n'est pas dans $L^2 (\R)$, donc les classes de cohomologie $L^2$ de $\alpha$ et $\beta$ sont non nulles. En revanche, si on utilise les coordonnées polaires $\xi=\rho\cos\theta$ et $\eta=\rho\sin\theta$, alors
\begin{eqnarray*}
\log(\frac{1}{|\xi|})=\log(\frac{1}{\rho})+\log(\frac{1}{|\cos\theta|})\geq \log(\frac{1}{\rho}),
\end{eqnarray*}
d'où
\begin{eqnarray*}
|\g_x|^2 &\leq& \frac{|\xi\a(\xi)\b(\eta)|^2}{(\xi^2 +\eta^2 )^2}\\
&\leq& |\log(\frac{1}{\rho\cos\theta})|^{-1}|\log(\frac{1}{\rho\sin\theta})|^{-1}\rho^{-2}\\
&\leq& |\log(\frac{1}{\rho})|^{-2} \rho^{-2},
\end{eqnarray*}
donc $\g_x$ est dans $L^2 (\R^2 ,d\xi\,d\eta)$, et il en est de même de $\g_y$. On conclut que $\gamma\in L^2$ et $d\gamma=\alpha\times\beta$, donc $\gamma\in\Omega^{2,2}(\R^2 )$, et la classe de $\alpha\times\beta$ est nulle.\qed

\bigskip

Pour remédier à cette difficulté, on introduit une condition sur une classe de torsion $L^p$, appelée \emph{robustesse}, qui garantit qu'elle reste non nulle après produit cartésien.

Soient $M_{1}$ et $M_{2}$ deux vari{\'e}t{\'e}s riemanniennes compl{\`e}tes. Supposons que $\tp{k_{1}}(M_{1})\not=0$ et
$\tp{k_{2}}(M_{2})\not=0$. D'apr{\`e}s le lemme \ref{stokes}, il existe des
formes ferm{\'e}es $L^{p}$ $e_{i}$ sur $M_{i}$ et des formes $L_{p'}$
$e'_{i,j}$ telles que
$$
\int_{M_{i}} e_{i}\wedge e'_{i}=1 \et
\n{de'_{i,j}}_{L^{p'}}\quad\hbox{tende vers}\quad 0.
$$
La forme ferm{\'e}e $L^{p}$ $\omega=e_{1}\wedge e_{2}$ sur $M_{1}\times M_{2}$ est elle non nulle en cohomologie ? Pour l'affirmer, il faudrait
contr{\^o}ler la norme de $d(e'_{1,j}\wedge e'_{2,j})$, c'est-{\`a}-dire non seulement celles de $de'_{1,j}$ et $de'_{2,j}$, mais aussi celles de $e'_{1,j}$ et de $e'_{2,j}$. On doit autoriser que $\n{e'_{1,j}}$ tende vers l'infini, mais moins vite que $\n{de'_{2,j}}$ ne tend vers 0. Ceci motive la d{\'e}finition suivante.

\begin{defi}
\label{to}
Soit $M$ une vari{\'e}t{\'e} riemannienne compl{\`e}te de dimension $n$. On note $\hop{k}(M)$ le sous-ensemble de $\hp{k}(M)$ form{\'e} des classes \emph{robustes}, i.e. qui contiennent une forme $\omega$ ayant la propri{\'e}t{\'e} suivante. Il existe une suite
$\omega'_{j}\in\opr{n-k}(M)$ telle que
\begin{enumerate}
\item les int{\'e}grales $\d \int_{M} \omega\wedge\omega'_{j}$ ne tendent
pas vers $0$ ;
\item les normes $\d\n{\omega'_{j}}_{L^{p'}}$ tendent vers $+\infty$
polynomialement en $j$ ;
\item les normes $\d\n{d\omega'_{j}}_{L^{p'}}$ tendent vers $0$
exponentiellement en $j$.
\end{enumerate}
Enfin, on note $\toop{k}(M)=\hop{k}(M)\cap\tp{k}(M)$.
\end{defi}

\begin{rem}
Cas où $p=2$.
\end{rem}
Dès que la torsion $L^2$ est non nulle, il y a des classes robustes. Il est possible que cela persiste pour tout $p$, mais je ne sais pas le montrer. J'en suis donc réduit à construire à la main ces classes robustes.

\begin{prop}
\label{direct}
Soient $M_{1}$ et $M_{2}$ des vari{\'e}t{\'e}s riemanniennes com\-pl{\`e}\-tes. Le produit cartésien de classes de torsion $L^p$ robustes de $M_1$ et $M_2$ respectivement est une classe robuste (et en particulier, non nulle) du produit riemannien $M_1 \times M_2$. Si l'une des deux classes est de degré maximum, le r{\'e}sultat est plus pr{\'e}cis : le produit cartésien est à nouveau une classe de torsion robuste.
\end{prop}

\preuve
Soient $\omega_{1}\in\op{k_1}(M_{1})$
(resp. $\omega_{2}\in\op{k_2}(M_{2})$) des formes ferm{\'e}es. Supposons
qu'il existe des formes $\omega'_{1,j}\in\opr{n_{1}-k_{1}}(M_{1})$
(resp. $\omega'_{2,j}\in\opr{n_{2}-k_{2}}(M_{2})$) comme dans la
d{\'e}finition \ref{to}. Posons $\omega=\pi_{1}^{*}\omega_{1}\wedge
\pi_{2}^{*}\omega_{2}$ et $\omega'_{j}=\pi_{1}^{*}\omega_{1,j}\wedge
\pi_{2}^{*}\omega_{2,j}$. Alors $\omega$ est ferm{\'e}e et $L^{p}$ et 
les formes $\omega'_{j}$ satisfont aux hypoth{\`e}ses de la d{\'e}finition
\ref{to}, donc la classe de cohomologie de $\omega$ est dans
$\hop{k_{1}+k_{2}}(M_{1}\times M_{2})$.

Supposons maintenant que $k_{1}=n_{1}$ et que 
$\omega_{2}$ est de torsion. Soit $\phi$ une $(n_{2}-k_{2})$-forme
ferm{\'e}e $L^{p'}$ sur $M_{1}\times M_{2}$. La restriction de $\phi$ {\`a}
presque tout facteur $\{*\}\times M_{2}$ est 
ferm{\'e}e et $L^{p'}$, donc pour presque tout $x_{1}\in M_{1}$,
$$
\int_{M_{2}}\omega_{2}\wedge\phi_{\mid\{x_{1}\}\times M_{2}}=0. 
$$
Il vient
$$
\int_{M_{1}\times M_{2}}\omega\wedge\phi
=\int_{M_{1}}\omega_{1}\wedge\int_{M_{2}}\omega_{2}\wedge\phi=0,
$$
autrement dit, $\omega$ est de torsion.\qed 

\subsection{Torsion des groupes ab{\'e}liens}

A titre d'application de la notion de classe robuste introduite au paragraphe précédent, montrons que la torsion $L^p$ de l'espace euclidien est non nulle en tout degré.

Commen{\c c}ons par le cas de la droite r{\'e}elle.

\begin{lemme}
\label{r}
L'ensemble $\toop{1}(\R)$ est non vide. Plus pr{\'e}cis{\'e}ment,
il existe une $1$-forme $L^{p}$ $a\,dt$ et une suite $u_j$ de fonctions
lisses {\`a} support compact sur $\R$ telles que
\begin{enumerate}
\item $\int_{\R} u_{j}a\,dt$ ne tend pas vers $0$ ;
\item $\n{u_{j}}_{L^{p'}}$ tend vers $+\infty$ polynomialement en $j$ ;
\item $\n{u'_{j}}_{L^{p'}}$ tend vers $0$ exponentiellement en $j$.

On a de plus les propri{\'e}t{\'e}s suivantes :

\item $\d\int_{\R} |sa'(s)|^{p}\,ds<+\infty$ ;
\item les fonctions $a$ et $-u_{j}$ sont d{\'e}croissantes sur $[1,+\infty[$
;
\item les fonctions $a$ et $u_{j}$ sont paires et s'annulent au
voisinage de $0$ ; 
\item $\n{u_{j}}_{L^{\infty}}$ tend vers $0$ exponentiellement
en $j$ ;
\item pour tout $\e>0$, $\n{s^{1-\e}u'_{j}}_{L^{p'}}$ et 
$\n{s^{-\e}u_{j}}_{L^{p'}}$ tendent vers $0$ exponentiellement.
\end{enumerate}
\end{lemme}

\preuve
Soit $\chi$ une fonction lisse et paire sur $\R$, {\`a} support dans
$[-1,1]$, qui vaut $1$
au voisinage de $0$.
On pose 
$$
a(x)=(1-\chi(x))|x|^{-{1\over{p}}}(\log |x|)^{-1} \si |x|\geq e,$$
$$
a(x)=(1-\chi(x))e^{-{1\over{p}}}\sinon.
$$
On d{\'e}finit une suite de fonctions $v_{j}$ paires, d{\'e}croissantes sur
$[0,+\infty[$ par
\begin{eqnarray*}
v_{j}(x)&=&2(1-\chi(x))2e^{-{j\over{p'}}} \si |x|\leq e^{j},\\
v_{j}(x)&=&2j\,|x|^{-{1\over{p'}}}(\log |x|)^{-1} \si e^{j}\leq |x|\leq
e^{2j},\\
v_{j}(x)&=&e^{-{{2j}\over{p'}}}(e^{j}+1-e^{-j}|x|)  \si e^{2j}\leq |x|\leq
e^{2j}(1+j^{-1}),\\
v_{j}(x)&=&0 \sinon.
\end{eqnarray*}
Comme $p>1$, $a$ et sa d{\'e}riv{\'e}e sont $L^p$. De plus $sa'(s)\sim
-{1\over{p}}a(s)$ donc $\int |sa'(s)|^{p}\,ds<+\infty$.
Par construction, $v'_{j}$
est nulle sur $[0,e^{j}]$, constante sur $[e^{j},e^{2j}(1+j^{-1})]$. Sur
l'intervalle $[e^{j},e^{2j}(1+j^{-1})]$, $\d |v'_{j}|$ est major{\'e}e par
$\d\con j\,s^{-1-1/p'}(\log s)^{-1}$. On calcule
$$
\int_{0}^{e^{j}} av_{j}= o(1),\quad
\int_{e^{j}}^{e^{2j}} av_{j}=1,\quad
\int_{e^{2j}}^{e^{2j}(1+j^{-1})} av_{j}=O(j^{-1}),
$$
$$
\int_{0}^{e^{j}} |v_{j}|^{p'} = 2^{p'},\quad
\int_{e^{j}}^{e^{2j}} |v_{j}|^{p'} =O(j),\quad
\int_{e^{2j}}^{e^{2j}(1+j^{-1})} |v_{j}|^{p'} =O(j^{-1}),
$$
$$
\int_{0}^{e^{j}} |v'_{j}|^{p'} = 0,\quad
\int_{e^{j}}^{e^{2j}} |v'_{j}|^{p'} =O(j^{p'}e^{-p'j}),
$$
$$
\int_{e^{2j}}^{e^{2j}(1+j^{-1})} |v'_{j'}|^{p'} =O(j^{p'-1}e^{-2p'j}),
\int_{1}^{e^{j}} |s^{-\e}v_{j}|^{p'} = O(je^{-\e p'j}),
$$
$$
\int_{e^{j}}^{e^{2j}} |s^{-\e}v_{j}|^{p'} = O(j^{p'}e^{-\e p'j}),\,
\int_{e^{2j}}^{e^{2j}(1+j^{-1})} |s^{-\e}v_{j}|^{p'} = O(j^{-1}e^{-2\e p'j}),
$$
$$
\int_{0}^{e^{j}} |s^{1-\e}v'_{j}|^{p'} = 0,\quad
\int_{e^{j}}^{e^{2j}} |s^{1-\e}v'_{j}|^{p'}=  O(j^{p'}e^{-\e p'j}) ,
$$
$$
\int_{e^{2j}}^{e^{2j}(1+j^{-1})} |s^{1-\e}v'_{j}|^{p'} = O(j^{p'-1}e^{-2\e p'j}).
$$
Une approximation $u_{j}$ lisse et {\`a} support compact de $v_{j}$ convient.\qed
 
\begin{cor}
\label{rn}
$\tp{k}(\R^{n})\not=0$ pour $k=1,\ldots,n$.
\end{cor}

\preuve
Montrons d'abord que $\toop1(\R^{n})$ est non vide.

Dans $\R^{n}$, on note $r$ la distance euclidienne {\`a}
l'origine. Soit $\theta'=dx_{2}\wedge\cdots\wedge dx_{n}$ une
$(n-1)$-forme parall{\`e}le. La forme $dr\wedge\theta$ {\'e}tant homog{\`e}ne de
degr{\'e} $0$, elle s'{\'e}crit
$$
dr\wedge\theta'=h\,dx_{1}\wedge\cdots\wedge dx_{n}
$$
o{\`u} la fonction $h$ est lisse en dehors de l'origine et homog{\`e}ne de degr{\'e}
$0$. Elle n'est pas identiquement nulle. Par homog{\'e}n{\'e}{\"\i}t{\'e}, $|h|$ et $r|dh|$ sont born{\'e}es.

Soient $a$ et $u_{j}$ les fonctions fournies par le
lemme \ref{r}. On consid{\`e}re les formes dif\-f{\'e}\-ren\-ti\-el\-les $\omega=a(r^{n})\,dr$
et $\omega'_{j}=u_{j}(r^{n})h\theta'$ sur $\R^{n}$. On v{\'e}rifie que
$$
\n{\omega}_{L^{p}}=\con(\int_{0}^{+\infty}
|a(r^{n})|^{p}r^{n-1}\,dr)^{1/p}=\con\n{a}_{L^{p}}
$$
est finie, que
$$
\int_{\R^{n}}\omega\wedge\omega'_{j}
=\int_{S^{n-1}}h^{2}\int_{0}^{+\infty} a(r^{n})u_{j}(r^{n})r^{n-1}\,dr
=C\,\int_{\R}au_{j}
$$
ne tend pas vers $0$, et que
$$
\n{\omega'_{j}}_{L^{p'}}=\con(\int_{0}^{+\infty}
|u_{j}(r^{n})|^{p'}r^{n-1}\,dr)^{1/p'}=\con\n{u_{j}}_{L^{p'}}
$$
tend vers $+\infty$ polynomialement.

On calcule
$$
d\omega'_{j}=nr^{n-1}u'_{j}(r^{n})h\, dr\wedge\theta'+
u_{j}(r^{n})\,dh\wedge\theta',
$$
et on majore
\begin{eqnarray*}
\n{r^{n-1}u'_{j}(r^{n})h}_{L^{p'}}
&\leq&\con(\int_{0}^{+\infty}
|r^{n-1}u'_{j}(r^{n})|^{p'}r^{n-1}\,dr)^{1/p'}\\
&=&\con\n{s^{n-1/n}u'_{j}(s)}_{L^{p'}},
\end{eqnarray*}
puis
\begin{eqnarray*}
\n{u_{j}(r^{n})\,dh}_{L^{p'}}
&\leq&\con(\int_{0}^{+\infty}
|r^{-1}u_{j}(r^{n})|^{p'}r^{n-1}\,dr)^{1/p'}\\
&=&\con\n{s^{-1/n}u_{j}(s)}_{L^{p'}},
\end{eqnarray*}
qui tendent vers $0$ exponentiellement, d'apr{\`e}s le lemme \ref{r}. On
conclut que $\omega$ est dans $\hop1(\R^{n})$, donc dans $\toop1(\R^{n})$
puisque la cohomologie r{\'e}duite est nulle.

Pour avoir le cas g{\'e}n{\'e}ral, il suffit d'appliquer
suffisamment de fois la proposition \ref{direct}.\qed

\subsection{Graduation des formes différentielles sur les produits semi-directs}

A la différence du cas des produits directs, une métrique riemannienne invariante à gauche sur un produit semi-direct $G=H\timez_{\alpha}\R$ croît exponentiellement, avec des exposants différents suivant les directions, déterminés par les valeurs propres de la dérivation $\alpha$. Ceci affecte les propriétés de contraction du champ de vecteurs invariant à droite $\xi$ qui engendre l'action à gauche du facteur $\R$ sur $G$. On a vu en section \ref{tor=0} comment utiliser le flot de ce champ de vecteurs - et ses propriétés de contraction de la norme $L^p$ de formes différentielles - pour montrer que la torsion $L^p$ est nulle. Les mêmes propriétés vont évidemment jouer un rôle dans la construction de classes de torsion non nulles.

\begin{defi}
\label{+}
Soit $G=H\timez_{\alpha}\R$ un produit semi-direct de groupes de Lie, soit $\H$ l'algèbre de Lie de $H$. Soit $p>1$ un réel, soit $k<\dim G$ un entier. \begin{enumerate}
  \item On décompose 
\begin{eqnarray*}
\Lambda^k \H^* =\Lambda^k_+ \oplus \Lambda^k_0 \oplus \Lambda^k_- ,
\end{eqnarray*}
où $\Lambda^k_+$ (resp. $\Lambda^k_0$, resp. $\Lambda^k_-$) est la somme des espaces caractéristiques de $\Lambda^k \alpha^\top$ relatifs aux valeurs propres de parties réelles supérieures (resp. égales, resp. inférieures) à $\frac{\tr\alpha}{p}$.
  \item On dit que $p$ est \emph{critique en degré $k$} pour $G$ si $\frac{\tr\alpha}{p}$ est la partie réelle d'une valeur propre de l'endomorphisme $\Lambda^k \alpha^\top$ de $\Lambda^k \H^*$, autrement dit, si $\Lambda^k_0 \not=0$.
  \item On note $d_0$, $d_{\pm}$ la différentielle extérieure composée avec le projecteur sur $\Lambda^k_0$, $\Lambda^k_{\pm}$.
\end{enumerate} 
\end{defi}
La décomposition dépendant de $p$, on notera $\Lambda^k_{+(p)}$, $d_{+(p)}$ s'il est nécessaire de spécifier l'exposant. 

\bigskip

\begin{rem}
Comparaison avec la définition \ref{kpcontractant}.
\end{rem}
Le champ $\xi$ est $(k,p)$-contractant si et seulement si $\Lambda^k_{+(p)}=\Lambda^k_{0(p)}=0$, $(k,p)$-dilatant si et seulement si $\Lambda^k_{-(p)}=\Lambda^k_{0(p)}=0$. Dans ce cas, d'après la section \ref{tor=0}, la torsion a des chances d'être nulle. Pour construire des classes de torsion, on va donc exploiter le fait que $\Lambda^k_{+(p)}$ et $\Lambda^k_{-(p)}$ sont simultanément non nuls.

\subsection{Critère de non-nullit{\'e} de la torsion}
\label{semi}

On s'inspire de la discussion des produits directs (paragraphe \ref{prod}).
Lorsqu'on passe aux produits semi-directs $G=H\timez_{\alpha}\R$ et qu'on s'in\-t{\'e}\-res\-se {\`a} un exposant $p$ non critique, on utilise seulement le fait que la cohomologie {\`a} support compact de la droite r{\'e}elle est non nulle. D'une certaine façon, l'op{\'e}rateur $\dplus$ remplace l'op{\'e}rateur $d$ sur l'autre facteur.

\begin{prop}
\label{tronc}
Soit $G=H\timez_{\alpha}\R$ un produit semi-direct. Soit $p>1$. Soit $e$ une $k-1$-forme ferm{\'e}e $L^{p}$ sur $H$. On suppose qu'il existe une suite de formes $e'_{j}\in\Omega^{n-k-1,p'}(H)$ telle que
\begin{enumerate}
\item les int{\'e}grales $\d\int_{H} e\wedge d_{+(p')} e'_{j}$ ne tendent pas
vers $0$ ;
\item la suite $\d m_{j}=\n{de'_{j}}_{L^{p'}}$ tend vers $+\infty$
polynomialement ;
\item la suite $\d n_{j}=\n{dd_{+(p')} e'_{j}}_{L^{p'}}$ tend vers $0$ exponentiellement.
\end{enumerate}
Alors $\hop{k}(G)\not=0$. Si de plus $H$ est nilpotent, alors $\toop{k}(G)\not=0$.
\end{prop}

\preuve
Notons $\pi:G\to H$ la projection dont les fibres sont les orbites de l'action à droite du facteur $\R$. Soit $\chi$ une fonction lisse sur $\R$ telle que $\chi=0$ au voisinage de $-\infty$ et $\chi=1$ au voisinage de $+\infty$. On note $\chi_{s}:t\mapsto \chi(t+s)$.

Posons
\begin{eqnarray*}
\omega=d\chi\wedge\pi^* e,
\end{eqnarray*}
La forme fermée $\omega$ représente le produit cartésien du générateur de la cohomologie à support compact $H_{c}^{1}(\R)$ et de la classe $[e]\in \tp{k-1}(H)$ (qui est non nulle, en vertu des hypothèses 1 et 3). Pour montrer que sa classe de cohomologie $L^p$ est non nulle, on utilise la dualité \ref{stokes} avec les formes-test $\omega'_{j}$ définies comme suit.
$$
\psi_{j}=\chi_{1}\pi^{*}\dplus e'_{j}+(1-\chi_{1})\pi^{*}\dmoins
e'_{j}+d\chi_{1}\wedge e'_{j}
$$
et
$$
\omega'_{j}=\chi_{s_{j}}(1-\chi_{-s_{j}})\psi_{j},
$$
o{\`u} $s_{j}$ est un r{\'e}el positif. Autrement dit, $\omega'_j$ est une troncature (destinée à rendre sa différentielle $L^{p'}$) d'une forme $\psi_j$ qui est $L^{p'}$, mais dont la différentielle ne l'est pas.

Alors
\begin{eqnarray*}
\int_{G}\omega\wedge\omega'_{j}
&=& \int_{G}d\chi\wedge e\wedge\psi_{j}\\
&=&\int_{\R}\chi_{1}d\chi \int_{H} e\wedge\dplus e'_{j}
-\int_{\R} (1-\chi_{1})d\chi \int_{H} e\wedge\dmoins e'_{j}\\
&=&\int_{H}e\wedge\dplus e'_{j}
\end{eqnarray*}
ne tend pas vers $0$. Comme
$\d d\psi_{j}=\pi^{*}d\dplus e'_{j}$, 
\begin{eqnarray*}
\n{\chi_{s_{j}}(1-\chi_{-s_{j}})d\psi_{j}}_{L^{p'}(G)}&\leq&
e^{\mu s_{j}}\n{d\dplus e'_{j}}_{L^{p'}(H)}\\
&\leq&e^{\mu s_{j}}n_{j}.
\end{eqnarray*}
D'autre part,
\begin{eqnarray*}
&&\n{d(\chi_{s_{j}}(1-\chi_{-s_{j}}))\wedge\psi_{j}}^{p'}_{L^{p'}(G)}\\
&=&
\n{d(\chi_{s_{j}})\wedge\psi_{j}}^{p'}_{L^{p'}(G)}
+\n{d(1-\chi_{-s_{j}})\wedge\psi_{j}}^{p'}_{L^{p'}(G)}\\
&\leq&\con(
\n{(\phi_{s_{j}})^{*}\dplus e'_{j}}^{p'}_{L^{p'}(H)}+
\n{(\phi_{-s_{j}})^{*}\dmoins e'_{j}}^{p'}_{L^{p'}(H)})\\
&\leq&\con
e^{-\eta p's_{j}}(\n{\dplus e'_{j}}^{p'}_{L^{p'}(H)}+
\n{\dmoins e'_{j}}^{p'}_{L^{p'}(H)})\\
&\leq&\con
e^{-\eta p's_{j}}\n{d e'_{j}}^{p'}_{L^{p'}(H)}.
\end{eqnarray*}
Il vient
$$
\n{\omega'_{j}}_{L^{p'}(G)}
\leq C\,(e^{\mu s_{j}}n_{j}+
e^{-\eta s_{j}}m_{j}).
$$
Posons $\d s_{j}=\inv{\mu+\eta}\log(\eta m_{j}/\mu n_{j})$. Avec ce choix, 
$$
\n{d\omega'_{j}}_{L^{p'}}\leq C'\,m_{j}^{\mu/\mu+\eta}n_{j}^{\eta/\mu+\eta},
$$
qui tend vers $0$ lorsque $j$ tend vers $+\infty$. Le lemme \ref{stokes}
entra{\^\i}ne alors que $\omega$ est non nulle dans $\hp{k}(G)$.

Soit $\omega'$ une $n-k$-forme ferm{\'e}e $L^{p'}$ sur $G$. Ecrivons
$\omega'=\beta'_{t}+dt\wedge \gamma'_{t}$. Alors $\beta'_{t}$ est une
forme ferm{\'e}e sur $H$ qui est dans $L^{p'}$ pour presque tout
$t\in\R$. Supposons $H$ unimodulaire. D'apr{\`e}s le th{\'e}or{\`e}me \ref{rpnul},
la cohomologie r{\'e}duite $\rp{k-1}(H)$ est nulle. Par cons{\'e}quent, pour
toute $k-1$-forme ferm{\'e}e $L^{p}$ $e$ sur $H$,
$\int_{H}e\wedge\beta'_{t}=0$. Il vient
$$
\int_{G} d\chi\wedge\pi^{*}e\wedge\omega'=\int_{\R}\chi'(t)\,dt\int_{H}
e\wedge\beta'_{t}=0.
$$
Ceci prouve que $\omega=d\chi\wedge\pi^{*}e$ est dans $\tp{k}(G)$, et
donc que $\tp{k}(G)\not=0$.\qed

\begin{rem}
Double valeur au bord.
\end{rem}
Notons $\xi$ le champ de vecteurs invariant à gauche qui engendre l'action à droite du facteur $\R$. Si on transporte $\psi_j$ par son flot, on trouve des limites distinctes, respectivement $d_{+}e'_{j}$ quand $t$ tend vers $+\infty$ et $d_{-}e'_{j}$ quand $t$ tend vers $-\infty$. Cela illustre le fait que, bien que $\xi$ ne soit ni $(k,p)$-contractant, ni $(k,p)$-dilatant, on peut définir deux valeurs au bord. Ce point de vue est développé dans \cite{PP}.

\subsection{Construction explicite de classes de cohomologie $L^p$}

Avant de se lancer dans le cas général, traitons un exemple.

\begin{exemple}
$T^{3,p}(G_{2,4,-\frac{1}{4}})\not=0$ si $\frac{4}{3}< p < 2$.
\end{exemple}
\preuve
Ici, $H={\mathbf{R}}^3$ avec les coordonnées $x,y,z$, et une matrice $\alpha$ diagonale, de valeurs propres 1, 1 et 2. Lorsque $\frac{4}{3}< p < 2$,  l'exposant conjugué $p'$ satisfait $2<p'<4$, d'où $2<\frac{\tr\alpha}{p'}<3$. Sur les 1-formes, $\Lambda^{1}_{+(p')}$ est engendré par $dz$, $\Lambda^{1}_{-(p')}$ par $dx$ et $dy$. On va appliquer le lemme \ref{tronc} avec des fonctions à support compact $e'_j$ qui ne dépendent que de la distance $r$ à l'origine, et une 2-forme fermée $e\in d\Lambda_-$, de la forme $e=d(f\beta)$ où $f$ est une fonction de $r$ et $\beta$ une 1-forme homogène de degré 1.
 
Comme $e=f'(r)\beta+f(r)d\beta$, et comme $d\beta$ est homogène de degré $0$, $|d\beta|$ est homogène de degré $-1$, d'où
\begin{eqnarray*}
\n{de}\lp &\leq&\con(\n{f'(r)}_{L^p (\R^3 )}+\n{r^{-1}f'(r)}_{L^p (\R^3 )} )\\
&\leq&\con((\int_{0}^{+\infty}|f'(r)|^p r^2 \,dr)^{1/p}+(\int_{0}^{+\infty}|r^{-1}f(r)|^p r^2 \,dr)^{1/p})).
\end{eqnarray*}

Notons $e'_j =w_j (r^2)$. Alors
\begin{eqnarray*}
de'_{j}=w'_j (r^2 )d(r^2 )=2w'_j (r^2 )(x\,dx+y\,dy+z\,dz),
\end{eqnarray*}
\begin{eqnarray*}
d_+ e'_j =2w'_j (r^2 )z\,dz,
\end{eqnarray*}
\begin{eqnarray*}
dd_+ e'_j =4w''_{j}(r^2 )(x\,dx+y\,dy+z\,dz)\wedge z\, dz,
\end{eqnarray*}
d'où
\begin{eqnarray*}
m_j =\n{de'_j}_{L^{p'}}=\con(\int_{0}^{+\infty}|rw'_j (r^2 )|^p r^2 \,dr)^{1/p},
\end{eqnarray*}
\begin{eqnarray*}
n_j =\n{dd_+ e'_j}_{L^{p'}}=\con(\int_{0}^{+\infty}|r^2 w''_j (r^2 )|^p r^2 \,dr)^{1/p},
\end{eqnarray*}

Le plus délicat à contrôler est l'intégrale $I_j =\int_{\R^3}e\wedge d_+ e'_j$. Comme $e'_j$ est à support compact, la formule de Stokes s'applique, et
\begin{eqnarray*}
I_j =\int_{\R^3}f\beta\wedge dd_+ e'_j .
\end{eqnarray*}
Choisir $\beta$ invariante par rotations (autour de l'origine, ou même seulement autour de l'axe $Oz$) est impossible, car $dr\wedge dd_+ e'_j =(xdx+ydy)\wedge dd_+ e'_j =0$. Il faut donc casser la symétrie, c'est pourquoi on choisit $\beta$ proportionnelle à $dy$. On calcule
\begin{eqnarray*}
dy\wedge dd_+ e'_j &=&-4w''_{j}(r^2 )xz\,dx\wedge dy\wedge dz\\
&=&r^2 w''(r^2 )h(x,y,z)\,dx\wedge dy\wedge dz,
\end{eqnarray*}
où $h(x,y,z)=\frac{-4xz}{r^2}$ est une fonction homogène de degré 0. C'est cette fonction qui entre comme ingrédient dans $e$ : on prend $\beta=h\,dy$, d'où $e=d(f(r)h\,dy)$.

Avec ce choix, il vient
\begin{eqnarray*}
\pm I_j &=&\int_{\R^3}e\wedge d_+ e'_j \\
&=&\int_{\R^3}f(r)h\,dy\wedge dd_+ e'_j \\
&=&-4\int_{\R^3}fh^2 r^2 w''_j (r^2 )\,dx\wedge dy\wedge dz\\
&=&-4(\int_{S^2}h^2 )\int_{0}^{+\infty}f(r)w''_j (r^2 )r^4 \,dr.
\end{eqnarray*}

Reste à trouver $f$ et $w_j$. Plutôt que de construire des fonctions adhoc, il suffit de prendre
\begin{eqnarray*}
f(r)=ra(r^3),\quad w_j (s)=-\int_{|s|}^{+\infty}t^{-1/2}u_j (t^{3/2})\,dt
\end{eqnarray*}
où $a$ et $u_j$ sont les fonctions obtenues au lemme \ref{r}. Les propriétés 4 à 8 de ce lemme sont là pour garantir que $n_j$ tend vers 0 exponentiellement, que $m_j$ tend vers l'infini au plus polynômialement et que$I_j$ ne tend pas vers 0 (voir ci-dessous). \qed

\bigskip

Passons au cas général (un peu plus général que le théorème B).

\begin{prop}
\label{tpknonnul}
On consid{\`e}re un produit semi-direct $G=H\timez_{\alpha}\R$ o{\`u}
$H=\R^{n-1}$ est ab{\'e}lien. On note $\l_{1}\leq\ldots\leq\l_{n-1}$ les parties r{\'e}elles des valeurs propres de $\alpha$. On note $w_{k}=\l_{1}+\cdots+\l_{k}$ et $W_{k}=\l_{n-1-k}+\cdots+\l_{n-1}$. Si $\d w_{k-1}<\frac{w_{n-1}}{p}<W_{k-1}$ et si $p$ est non critique en degr{\'e} $k-1$, alors $\tp{k}(G)\not=0$.
\end{prop}

\preuve
Les in{\'e}galit{\'e}s $w_{k-1}<w_{n-1}/p<W_{k-1}$ entra{\^\i}nent que $w_{n-k}$ $<w_{n-1}/p'<W_{n-k}$.
Etant donn{\'e} $I\subset\{1,\ldots,n-1\}$, on note $\l_{I}=\sum_{i\in I}
\l_{i}$. Consid{\'e}rons, parmi les parties $I$ {\`a} $n-k$ {\'e}l{\'e}ments de
$\{1,\ldots,n-1\}$ telles que $\l_{I}>w_{n-1}/p'$, celle, not{\'e}e $I_{0}$, pour laquelle
$\l_{I}$ est minimum. Notons $i_{m}$ le plus petit indice qui n'est pas
dans $I_{0}$ et $i_{M}$ le plus grand {\'e}l{\'e}ment de $I_{0}$.
Comme $\l_{I_{0}}>w_{n-1}/p'>w_{n-k}$, $i_{m}\leq n-k$ et
$\l_{i_{m}}<\l_{i_{M}}$. Posons
$I_{1}=(I_{0}\cup\{i_{m}\})\setminus\{i_{M}\}$. Alors
$\l_{I_{1}}<\l_{I_{0}}$ donc par d{\'e}finition de $I_{0}$, $\l_{I_{1}}\leq
w_{n-1}/p'$. Comme $p$ est non critique en degr{\'e} $k-1$, $p'$ est non
critique en degr{\'e} $n-k$, donc $w_{n-1}-p'\l_{I}\not=0$ pour tout ensemble
$I$ {\`a} $n-k$ {\'e}l{\'e}ments. Par cons{\'e}quent, $\l_{I_{1}}<w_{n-1}/p'$.

Soit $\theta'\in\L^{n-k-1}\H^{*}$ un
vecteur propre de $\L^{n-k-1}\alpha$ relatif {\`a} une valeur propre de
partie r{\'e}elle $\mu'=\l_{I_{0}}-\l_{i_{M}}$, et soient $\eta$ et $\eta'\in
H^{*}$ des vecteurs propres relatifs {\`a} des valeurs propres de parties
r{\'e}elles $\l_{i_{m}}$ et $\l_{i_{M}}$ respectivement. Alors
$(\eta'\wedge\theta')_{+(p')}=\eta'\wedge\theta'$ mais
$(\eta\wedge\theta')_{+(p')}=0$ donc
$$
(\eta+\eta')\wedge((\eta+\eta')\wedge\theta')_{+(p')}=\eta\wedge\eta'\wedge\theta'
$$
est non nul. Il existe donc $\theta\in\L^{k-2}\H^{*}$ tel que 
$\d\theta\wedge(\eta+\eta')\wedge((\eta+\eta')\wedge\theta')_{+(p')}\not=0$.

Soient $a$ et $u_{j}$ les fonctions fournies par le
lemme \ref{r}. Posons, pour $s>0$, 
$$ 
{\tilde v}_{j}(s)=s^{-1/2}u_{j}(s^{(n-1)/2}).
$$
Notons $w_{j}$ la fonction {\`a} support compact sur $[0,+\infty[$ dont la d{\'e}riv{\'e}e est ${\tilde v}_{j}$. Dans $H=\R^{n-1}$, on note $r$ la distance euclidienne {\`a}
l'origine. On d{\'e}finit des fonctions $f=ra(r^{n-1})$ et $g_{j}=w_{j}(r^{2})$ sur $H$. Par construction, $dg_{j}=2u_{j}(r^{n-1})\,dr$. 

On consid{\`e}re les formes diff{\'e}rentielles $e'_{j}=g_{j}\theta'$ sur
$H$.
On a $de'_{j}=w'_{j}(r^{2})d(r^{2})\wedge\theta'$ donc $\dplus
e'_{j}=w'_{j}(r^{2}) dr_{+}^{2}\wedge\theta'$ o{\`u} on a not{\'e} 
$$
r_{+}^{2}=\sum_{\l_{i}+\mu'>\tr\alpha/p'} x_{i}^{2}.
$$
Comme la forme $d(r_{+}^{2})$ est ferm{\'e}e, 
$$
d\dplus e'_{j}=w''_{j}(r^{2})d(r^{2})\wedge d(r_{+}^{2})\wedge\theta'.
$$
Comme la $n-1$-forme $d(r^{2})\wedge\dplus
r^{2}\wedge\theta\wedge\theta'$ est homog{\`e}ne, on peut l'{\'e}crire
$$
d(r^{2})\wedge\dplus
r^{2}\wedge\theta\wedge\theta'
=r^{2}h(x)dx_{1}\wedge\cdots\wedge dx_{n-1}
$$
o{\`u} la fonction $h$ est lisse en dehors de l'origine et homog{\`e}ne de degr{\'e}
$0$. Par homog{\'e}n{\'e}{\"\i}t{\'e}, $|h|$ et $r|dh|$ sont born{\'e}es.

Il existe un point de $H$ o{\`u} $d(r^{2})=\eta+\eta'$. En ce point, la
$n-1$-forme $\theta\wedge d(r^{2})\wedge(d(r^{2})\wedge\theta')_{+(p')}$
est non nulle, donc $h$ n'est pas identiquement
nulle. On pose 
$$
e=d(fh\theta).
$$
Comme $e'_{j}$ est {\`a} support compact, 
\begin{eqnarray*}
\pm\int_{H}e\wedge \dplus e'_{j}&=&
\int_{H}fh\theta\wedge d\dplus e'_{j}\\
&=&
\int_{H} fw''_{j}(r^{2})h\,d(r^{2})\wedge\dplus
r^{2}\wedge\theta\wedge\theta'\\
&=&
\int_{S^{n-2}} h^{2}\int_{0}^{+\infty}
f(r)w''_{j}(r^{2})r^{2}r^{n-2}\,dr\\
&=& C\,\int_{0}^{+\infty} f(r)w''_{j}(r^{2})r^{n}\,dr
\end{eqnarray*}
o{\`u} $C>0$. On calcule, pour $r>0$,
$$
w''_{j}(r)=-{1\over2}r^{-3/2}u_{j}(r^{n-1/2})+{{n-1}\over2}r^{n-4/2}u'_{j}(r^{n-1/2}).
$$
Comme $u_{j}$ est d{\'e}croissante sur $[1,\infty[$, les deux termes de la
somme sont de m{\^e}me signe, donc
\begin{eqnarray*}
|\int_{1}^{+\infty} f(r)w''_{j}(r^{2})r^{n}\,dr|
&\geq&
{1\over2}\int_{1}^{+\infty} ra(r^{n-1})r^{-3}u_{j}(r^{n-1})r^{n}\,dr\\
&=&{1\over2}\int_{1}^{+\infty} a(s)u_{j}(s)\,ds\\
&=&{1\over4}\int_{\R} a(s)u_{j}(s)\,ds
\end{eqnarray*}
qui ne tend pas vers $0$. L'int{\'e}grale $\d \int_{0}^{1}
f(r)w''_{j}(r^{2})r^{n}\,dr$ tend vers $0$, donc 
$\d \int_{H}e\wedge \dplus e'_{j}$ ne tend pas vers $0$.

Comme $e=hf'(r)dr\wedge\theta+fdh\wedge\theta$,
\begin{eqnarray*}
\n{e}_{L^{p}}&\leq&\con \n{f'(r)}_{L^{p}(H)}+\n{r^{-1}f(r)}_{L^{p}(H)}\\
&=&
\con (\int_{0}^{+\infty} |f'(r)|^{p}r^{n-2}\,dr)^{1/p}\\
&&+\int_{0}^{+\infty}
|r^{-1}f(r)|^{p}r^{n-2}\,dr)^{1/p}\\
&\leq&
\con (\int_{0}^{+\infty} |r^{n-1}a'(r^{n-1})|^{p}r^{n-2}\,dr)^{1/p}\\
&&+\int_{0}^{+\infty}
|a(r^{n-1})|^{p}r^{n-2}\,dr)^{1/p}\\
&=&
\con (\int_{0}^{+\infty} |sa'(s)|^{p}\,ds)^{1/p}+\int_{0}^{+\infty}
|a(s)|^{p}\,ds)^{1/p}\\
\end{eqnarray*}
est finie.

Comme $de'_{j}=2u_{j}(r^{n-1})dr\wedge\theta'$,
\begin{eqnarray*}
\n{de'_{j}}^{p'}_{L^{p'}}
&=&\con\int_{0}^{+\infty} |u_{j}(r^{n-1})|^{p'}r^{n-1}\,dr\\
&=&\con\int_{0}^{+\infty} |u_{j}(s)|^{p'}\,ds
\end{eqnarray*}
cro{\^\i}t polynomialement en $j$, d'apr{\`e}s le lemme \ref{r}.

Comme $\d d\dplus e'_{j}=w''_{j}(r^{2})d(r^{2})\wedge\dplus
r^{2}\wedge\theta'$,
\begin{eqnarray*}
\n{d\dplus e'_{j}}_{L^{p'}}& \leq& \con
((\int_{0}^{+\infty} |r^{-1}u_{j}(r^{n-1})|^{p'}r_{n-1}\,dr)^{1/p'}\\
&&+
(\int_{0}^{+\infty}
|r^{n-2}u'_{j}(r^{n-1})|^{p'}r_{n-1}\,dr)^{1/p'})\\
&=&\con
((\int_{0}^{+\infty} |s^{-1/n-1}u_{j}(s)|^{p'}ds)^{1/p'}\\
&&+
(\int_{0}^{+\infty}
|s^{n-2/n-1}u'_{j}(s)|^{p'}ds)^{1/p'}),\\
\end{eqnarray*}
tend vers $0$ exponentiellement, d'apr{\`e}s le lemme \ref{r}. De la proposition \ref{tronc}, il r{\'e}sulte que $\tp{k}(G)\not=0$.\qed

\subsection{Torsion de l'espace hyperbolique réel}
\label{torhyp}

Lorsqu'on s'in\-t{\'e}\-res\-se {\`a} un exposant $p$ critique, on utilise le fait que la cohomologie $L^p$ de la droite r{\'e}elle est non nulle, ainsi qu'une information plus fine sur la cohomologie $L^p$ de $H$, faisant jouer un rôle important à l'op{\'e}rateur $\dzero$.

\begin{prop}
\label{torcri}
Soit $G=H\timez_{\alpha}\R$ un produit semi-direct. Soit $p$ un exposant
critique en degr{\'e} $k-1$, i.e. tel que $\L^{k-1}_{0}\not=0$. On suppose
qu'il existe une $k-1$-forme $e$ ferm{\'e}e et $L^{p}$ sur $H$ et une suite 
$e'_{j}\in\Omega^{n-k,p'}(H)$ telles que
\begin{itemize}
\item $e_{-}=0$, $e'_{j,-}=0$, $\dmoins e'_{j}=0$ ;
\item $\int_{H} e\wedge e'_{j}$ ne tend pas vers $0$ ;
\item $\n{e'_{j,0}}_{L^{p'}}$, $\n{e'_{j,+}}_{L^{p'}}$ et $\n{\dplus e'_{j}}_{L^{p'}}$ tendent vers l'infini poly\-n{\^o}\-mia\-le\-ment en $j$ ;
\item $\n{\dzero e'_{j}}_{L^{p'}}$ tend vers $0$ exponentiellement en $j$.
\end{itemize}
Alors $\hp{k}(G)\not=0$.
\end{prop}

\preuve
Soient $a$ et $u_{j}$ des fonctions sur $\R$ qui sont nulles sur
$[0,+\infty[$ et co{\"\i}ncident sur $]-\infty,0]$ avec celles construites en
\ref{r}. 

On pose $\omega=a(t)\,dt\wedge\pi^{*}e$ et
$\omega'_{j}=u_{j}(t)\pi^{*}e'_{j}$. Alors $\omega$ est ferm{\'e}e. Comme
$e_{-}=0$, il
existe une constante strictement positive $\nu$ telle que
$$
\n{\omega}_{L^{p}(G)}^{p}
\leq\con 
\int_{-\infty}^{0} |a(t)|^{p}(
e^{\nu t}\n{e_{+}}_{L^{p}(H)}^{p}+\n{e_{0}}_{L^{p}(H)}^{p})\,dt,
$$
donc $\omega\in L^{p}(G)$. De m{\^e}me,
$$
\n{\omega'_{j}}_{L^{p'}(G)}^{p'}
\leq\con 
\int_{-\infty}^{0} |u_{j}(t)|^{p'}(
e^{\nu t}\n{e'_{j,+}}_{L^{p'}(H)}^{p'}+\n{e_{j,0}}_{L^{p'}(H)}^{p'})\,dt
$$
tend vers $+\infty$ polyn{\^o}mialement. On calcule
$$
d\omega'_{j}=u'_{j}(t)\,dt\wedge\pi^{*}e'_{j}+u_{j}(t)\pi^{*}de'_{j},
$$
et 
\begin{eqnarray*}
&&\n{u'_{j}(t)\,dt\wedge\pi^{*}e'_{j}}_{L^{p'}(G)}^{p'}\\
&\leq&\con 
\int_{-\infty}^{0} |u'_{j}(t)|^{p'}(
e^{\nu t}\n{e'_{j,+}}_{L^{p'}(H)}^{p'}+\n{e_{j,0}}_{L^{p'}(H)}^{p'})\,dt
\end{eqnarray*}
tend vers $0$ exponentiellement car $\n{u'_{j}}_{L^{p'}}$ et
$\n{u'_{j}}_{L^{\infty}}$ tendent vers $0$ exponentiellement. De m{\^e}me
\begin{eqnarray*}
&&\n{u_{j}(t)\pi^{*}de'_{j}}_{L^{p'}(G)}^{p'}\\
&\leq&\con 
\int_{-\infty}^{0} |u_{j}(t)|^{p'}(
e^{\nu t}\n{\dplus e'_{j}}_{L^{p'}(H)}^{p'}+\n{\dzero e_{j}}_{L^{p'}(H)}^{p'})\,dt
\end{eqnarray*}
tend vers $0$ exponentiellement car $\n{\dzero e'_{j}}_{L^{p'}}$ et
$\n{u_{j}}_{L^{\infty}}$ tendent vers $0$ exponentiellement.

Enfin 
$$
\int_{G}\omega\wedge\omega'_{j}
=\int_{-\infty}^{0} au_{j}\int_{H}e\wedge e'_{j}
$$
ne tend pas vers $0$. On conclut avec le lemme \ref{stokes} que $\omega$
est non nulle dans $\hp{k}(G)$.\qed

\begin{cor}
\label{hyperbreeltor}
Soit $M=\R H^n$ l'espace hyperbolique r{\'e}el de dimension $n$. Pour chaque $2\leq k\leq n-1$, 
$$
\tp{k}(M)\not=0\quad\Leftrightarrow\quad p=\frac{n-1}{k-1}.
$$
\end{cor}

\preuve
Le théorème A s'applique et entraîne que, pour tout $p\not=\frac{n-1}{k-1}$, $\tp{k}(M)=0$. De plus, $\hp{k}(M)=0$ dès que $p\leq\frac{n-1}{k-1}$.

Réciproquement, soit $e$ une $k-1$-forme ferm{\'e}e {\`a} support compact
sur $\R^{n-1}$, non nulle. Cela existe d{\`e}s que $k\geq 2$. Soit
$e'_{j}=e'$ une $n-k$-forme sur $\R^{n-1}$ telle que $\int e\wedge
e'\not=0$. Comme $\L^{k-1}=\L^{k-1}_{0}$, $\L^{n-k}=\L^{n-k}_{0}$ et
$\L^{n-k+1}=\L^{n-k+1}_{+}$, les conditions $e_{-}=0$, $e'_{-}=0$ et $\dmoins e'=0$ sont automatiquement satisfaites. De plus, $\dplus e'=0$, donc la proposition \ref{torcri} s'applique, et $\hp{k}(M)\not=0$. 

Pour montrer que $\rp{k}(M)=0$, on utilise la dualité de Poincaré, corollaire \ref{poincare}. En degré $k'=n-k$, l'exposant conjugu{\'e} $\d p'=\frac{n-1}{n-k}=\mathbf{q}{n,-1,k'}$, est justement le cas limite d'application du théorème A, donc $H^{n-k,p'}(M)=0$. En particulier, $R^{n-k,p'}(M)=0$, d'où, par dualité, $\rp{k}(M)=0$. On conclut que $\tp{k}(M)\not=0$.\qed

\subsection{Preuve du théorème B}

Le cas de l'espace hyperbolique réel a fait l'obet du paragraphe \ref{torhyp}.

Le groupe $G_{\mu,n,\delta}$ s'obtient en faisant $\l_1 =\cdots=\l_{n-\mu}=\sqrt{-\delta}$, $\l_{n-\mu+1} =\cdots=\l_{n-1}=1$. On pose $k=\mu$, il vient $W_{k-1}=k-1$, $w_{n-1}=k-1+(n-k)\sqrt{-\delta}$. Les exposants critiques en degré $k-1$ sont les nombres de la forme $w_{n-1}/\lambda$ où $\lambda$ est une somme de $k-1$ nombres parmi $\l_1 ,\ldots,\l_{n-1}$. Le plus petit est $\d\frac{w_{n-1}}{W_{k-1}}$. Le suivant est $\d\frac{w_{n-1}}{W_{k-2}+\l_{n-k}}>\frac{w_{n-1}}{w_{k-1}}$. On peut donc appliquer la proposition \ref{tpknonnul} et conclure que $\tp{k}(G)\not=0$ pour tout $p$ dans l'intervalle 
\begin{eqnarray*}
]\frac{w_{n-1}}{W_{k-1}},\frac{w_{n-1}}{W_{k-2}+\l_{n-k}}[\quad =\quad ]{\bf q}(n,\delta,k-1),1+\frac{1+(n-k-1)\sqrt{-\delta}}{k-2+\sqrt{-\delta}}[.
\end{eqnarray*}

Si $G_{\mu,n,\delta}$ était quasiisométrique à une variété riemannienne $M$ simplement connexe, complète, à courbure sectionnelle négative $\delta'$-pincée pour un $\delta'<\delta$ proche de $\delta$, alors, pour ${\bf q}(n,\delta,\mu-1)<p<{\bf q}(n,\delta',\mu-1)$, $\tp{\mu}(G_{\mu,n,\delta})\not=0$ mais le théorème A donne que $\tp{\mu}(M)=0$, c'est incompatible avec l'invariance sous quasiisométrie de la cohomologie $L^p$ pour les espaces uniformément contractiles, voir \cite{Gromovasympt}, section 8.\qed

\begin{exemple}
Cas des espaces symétriques de rang un.
\end{exemple}
La construction qui précède ne s'étend pas aux produits semi-directs $G=H\timez_{\alpha}\R$ où $H$ est nilpotent non abélien. On l'explique sur l'exemple où $G$ est isométrique au plan $H^2_{\mathbf{C}}$. Dans ce cas, $H$ est le groupe d'Heisenberg de dimension 3. Son algèbre de Lie $\H$ admet une base $(X,Y,Z)$ où $Z=[X,Y]$ est central. Dans cette base, la matrice de la dérivation $\alpha$ est 
$\d\begin{pmatrix}1 & 0 & 0 \\ 0 & 1 & 0\\ 0 & 0 & 2\end{pmatrix}$. 
Soient $dx$, $dy$, $\tau = dz-x dy$ les éléments de la base duale, vus comme formes invariantes à gauche sur $H$. Ce sont des vecteurs propres de $\alpha$, pour les valeurs propres $-1$, $-1$ et $-2$. Soit $p$ un réel, $\frac{4}{3}<p<2$. On s'intéresse à la torsion en degré 2, $\tp{2}(G)$. Alors $2<p'<4$,  $1<\frac{\tr \alpha}{p'}<2$, donc $\Lambda_{+(p')}$ est engendré par $\tau$, $\Lambda_{-(p')}$ par $dx$ et $dy$.

Soit $e'_j$ une suite de fonctions sur $H$. Alors
\begin{eqnarray*}
de'_j & = & (Xe'_j) dx + (Ye'_j) dy + (Ze'_j) \tau \\
d_+ e'_j & = & (Ze'_j)\tau  \\
dd_+ e'_j & = & d(Ze'_j) \wedge \tau  + (Ze'_j)d\tau \, .
\end{eqnarray*}
En particulier, $Ze'_j = - dd_+ e'_j(X\wedge Y)$.

Par conséquent,
$$
|\n{d_+ e'_j} _{L^{p'}} \leq  \n{dd_+ e'_j}_{L^{p'}},
$$
et $\int _{H} e \wedge d_+ e'_j$ doit tendre vers 0 en même temps que $\n{dd_+ f_j}_{L^{p'}}$. Il n'existe donc aucune donnée $e$, $e'_j$ qui satisfasse aux hypothèses du lemme \ref{tronc}. 

De fait, $\tp{2}(\mathbf{C}H^2)=0$ pour $\frac{4}{3}<p<2$, voir \cite{PP}.

\vskip1cm
\noindent Laboratoire de Math{\'e}matique d'Orsay\\
UMR 8628 du C.N.R.S.\\
Universit{\'e} Paris-Sud 11\\
B{\^a}timent 425\\
91405 Orsay\\
France\\
\smallskip
{\tt\small Pierre.Pansu@math.u-psud.fr\\
http://www.math.u-psud.fr/$\sim$pansu}


\begin{thebibliography}{}

\bibitem[Be]{Berger} M. BERGER,
{\sl Sur certaines vari\'et\'es riemanniennes \`a courbure positive.}
C. R. Acad. Sci., Paris ${\bf 247}$, $1165-1168 ~ (1958)$.

\bibitem[BK]{courbure} P. BUSER and H. KARCHER,
{\sl Gromov's almost flat manifolds.}
Ast{\'e}risque ${\bf 81}$, Soc. Math. de France, Paris $(1981)$.

\bibitem[CE]{CE} J. CHEEGER, D. EBIN,
{\sl Comparison theorems in Riemannian geometry.}
Amsterdam: North Holland $(1975)$.

\bibitem[DX]{Donnelly-Xavier} H. DONNELLY, F. XAVIER,
{\sl On the differential form spectrum for negatively curved manifolds.}
Amer. J. Math. $\bf 108$, $169-185 ~ (1984)$.

\bibitem[G2]{Gromovasympt} M. GROMOV,
{\sl Asymptotic invariants of infinite groups.}
In ``Geometric Group Theory'', ed. G. Niblo and M. Roller, Cambridge: Cambridge University Press, $(1993)$.

\bibitem[GKS]{GKS} V. GOLDSTEIN, V. KUZMINOV, I. SHVEDOV,
{\sl The Kuenneth formula for $L_{p}$ cohomologies of warped products.}
Sib. Math. J. ${\bf 32}$, No.5, $749-760 ~ (1991)$; translation from
Sib. Mat. Zh. ${\bf 32}$, No.5(189), $29-42 ~ (1991)$.

\bibitem[GT]{GT} V. GOLDSTEIN, M. TROYANOV, 
{\sl The $L^{p,q}$ cohomology of $SOL$.}
Ann. Fac. Sci. Toulouse ${\bf 7}$,  $687-689 ~ (1998)$. 

\bibitem[He]{Heintze} E. HEINTZE,
{\sl On homogeneous manifolds of negative curvature.}
Math. Annalen ${\bf 211}$, $23-34 ~ (1974).$

\bibitem[Hz]{Hernandez} L. HERN{\'A}NDEZ LAMONEDA,
{\sl K{\"a}hler manifolds and $\inv4$-pinching.}
Duke Math. J. ${\bf 62}$, $601-611 ~ (1991)$.

\bibitem[JY]{JY} J. JOST, S.T. YAU,
{\sl Harmonic maps and superrigidity.}
Greene, Robert (ed.) et al., Differential geometry. Part 1: Partial differential equations on manifolds. Proceedings of a
summer research institute, held at the University of California, Los Angeles, CA, USA, July 8-28, 1990. Proc. Symp. Pure Math. ${\bf 54}$, Part 1, $245-280 ~ (1993)$.

\bibitem[KS]{KS} V. KUZMINOV, I. SHVEDOV, 
{\sl On compact solvability of the operator of exterior derivation.}
Sib. Math. J. ${\bf 38}$, No.3, $492-506 ~ (1997)$; translation from
Sib. Mat. Zh. ${\bf 38}$, No.3, $573-590 ~ (1997)$.

\bibitem[MSY]{MSY} N. MOK, Y.T. SIU, S.K. YEUNG,
{\sl Geometric superrigidity.} 
Invent. Math. ${\bf 113}$, $57-83 ~ (1993)$.

\bibitem[P1]{P'} P. PANSU,
{\sl Differential forms and connections adapted to contact structures,
after M. Rumin.}
P. 183-196 in ``Symplectic Geometry'', D. Salamon ed., L.M.S. Lect. Notes
Vol. ${\bf 192}$, London: London Math. Soc. $(1993)$. 

\bibitem[P2]{PP} P. PANSU,
{\sl Cohomologie $L^p$, espaces homog{\`e}nes et pincement.}
Texte disponible à la page \texttt{http://www.math.u-psud.fr/\%7Epansu/liste-prepub.html} depuis $(1999)$.

\bibitem[P3]{PCambridge} P. PANSU,
{\sl $L^p$-cohomology and pinching.}
M. Burger, A. Iozzi (ed.), Rigidity in Dynamics and Geometry. Proc. Cambridge 2000. Springer Verlag, Berlin... $379-389 ~ (2002)$.

\bibitem[V]{Ville} M. VILLE,
{\sl On $\frac{1}{4}$-pinched 4-dimensional Riemannian manifolds of negative curvature.}
Ann. Global Anal. Geom. ${\bf 3}$, $329-336 ~ (1985)$.

\bibitem[YZ]{YZ} S.T. YAU and F. ZHENG,
{\sl Negatively $\frac{1}{4}$-pinched Riemannian metric on a compact Kähler manifold.}
Invent. Math. ${\bf 103}$ $527-535 ~ (1991)$.

\end{thebibliography}
\end{document}